\documentclass[preprint,12pt,times]{elsarticle}
\usepackage{lineno}
\modulolinenumbers[200]
\usepackage{color}
\usepackage{colortbl}
\usepackage{bm}
\usepackage{transparent}
\usepackage{graphicx}
\usepackage{import}
\usepackage{diagbox}
\usepackage{algorithmic}
\usepackage{algorithm}
\usepackage{appendix}
\usepackage{amsmath}
\usepackage{amscd}
\usepackage{mathrsfs}
\usepackage{graphicx,color}
\usepackage{arydshln}
\usepackage{booktabs}
\usepackage{multirow}
\usepackage{rotating}
\usepackage{setspace}
\usepackage{indentfirst}
\usepackage{cases}
\usepackage{verbatim}
\usepackage{xcolor}
\usepackage{float}

\graphicspath{{./figure/}}
\usepackage{amssymb,color}
\usepackage{amsthm}
\usepackage{epsfig,subfigure}

\newcommand{\edge}[1]{\ar@{-}[#1]}
\numberwithin{equation}{section}

\biboptions{numbers,sort&compress}%
\begin{document}

\begin{frontmatter}
    %\title{Data-Free Deep Learning Approach for Solving the Inverse Problem of the Wave Equation}
    \title{ An Unsupervised Deep Learning Approach for the Wave Equation Inverse Problem}

    \author[1,4]{Xiong-Bin Yan}
    \ead{yanxb2015@163.com}
    \author[1]{Keke Wu}
    \author[1,2,3]{Zhi-Qin John Xu\corref{cor}}
    \ead{xuzhiqin@sjtu.edu.cn}
    \author[1,2,3,4]{Zheng Ma\corref{cor}}
    \ead{zhengma@sjtu.edu.cn}
    \cortext[cor]{Corresponding author}
    \address[1]{School of Mathematical Sciences, Shanghai Jiao Tong University, Shanghai, China}
    \address[2]{Institute of Natural Sciences, MOE-LSC, Shanghai Jiao Tong University, Shanghai, China}
    \address[3]{Qing Yuan Research Institute, Shanghai Jiao Tong University, Shanghai, China}
    \address[4]{CMA-Shanghai, Shanghai Jiao Tong University, Shanghai, China}
    \begin{abstract}
        Full-waveform inversion (FWI) is a powerful geophysical imaging technique that infers high-resolution subsurface physical
        parameters by solving a non-convex optimization problem. However, due to limitations in observation, e.g., limited
        shots or receivers, and random noise, conventional inversion methods are confronted with numerous
        challenges, such as the local-minimum problem. In recent years, a substantial body of work has demonstrated that
        the integration of deep neural networks and partial differential equations for solving full-waveform inversion problems
        has shown promising performance. In this work, drawing inspiration from the expressive capacity of neural networks,
        we provide an unsupervised learning approach aimed at accurately reconstructing subsurface physical
        velocity parameters. This method is founded on a re-parametrization technique for Bayesian inference,
        achieved through a deep neural network with random weights. Notably, our proposed approach does not hinge
        upon the requirement of the labeled training dataset, rendering it exceedingly versatile and adaptable
        to diverse subsurface models. Extensive experiments show that the proposed approach performs noticeably better
        than existing conventional inversion methods.
    \end{abstract}
    \begin{keyword}
Full-waveform inversion; Bayesian inference; Unsupervised learning; Re-parametrization technique; Deep neural network.
    \end{keyword}
\end{frontmatter}

\section{Introduction}\label{Sec1}
In geophysics, seismic waveform inversion is frequently employed to obtain quantitative
estimates of subsurface properties that can accurately predict the observed seismic data.
The reconstruction of subsurface properties is a profoundly non-linear inverse problem.
There are various techniques available to address the seismic inverse problem, such as
velocity analysis of stacked seismic traces \cite{Berkhout+1997}, migration-based travel-time
approaches \cite{zelt+1992,clement+2001}, Born approximation \cite{hudson1981use, muhumuza2018seismic},
and full-waveform inversion \cite{tarantola1984inversion, warner2013anisotropic, jakobsen2015full}.
Distinguished from other approaches, full-waveform inversion excels in the ability
to deduce high-resolution subsurface structures. This is achieved through an iterative process
that involves aligning observed and simulated seismograms while harnessing comprehensive wavefield
information. This superiority stems from its capacity to fully exploit the informational richness
of recorded seismic data, encompassing both amplitude and travel-time components.

In the context of full-waveform inversion, the typical approach for
minimizing the objective function involves the use of gradient descent methods.
These methods require the explicit computation of gradients pertaining to the cost function
with respect to the velocity model. Nevertheless, the application of gradient descent methods
to the full-waveform inversion (FWI) objective function poses considerable challenges. This
is primarily due to the presence of numerous local minima, a consequence of the
inherent high nonlinearity and ill-posed nature of the problem \cite{virieux2009overview}.
Over the past few decades, researchers have developed numerous approaches aimed at enhancing
the effectiveness of FWI. These approaches encompass the integration of multiple data
components \cite{bunks1995multiscale}, the introduction of regularization terms \cite{asnaashari2013regularized,yang2023wasserstein},
and the adoption of novel objective functions \cite{yang2018application,qiu2017full,li2022quadratic}.

Deep learning, primarily in the form of Deep Neural Networks (DNNs), has
gained significant attention in the fields of science and engineering. Its versatile applications span a wide range,
including image classification/recognition \cite{zou2015deep}, shape representation \cite{wu20153d},
and natural language processing \cite{deng2018deep}. Similarly, this method has also
found extensive application in the solution of inverse problems and has demonstrated
performance surpassing that of traditional methods. Examples include CT reconstruction \cite{wang2020deep,ding2022dataset},
image processing \cite{chen2022nonblind, tian2020deep}, MRI reconstruction \cite{hyun2018deep}, and more.

A strong effort has been made in recent years to facilitate deep learning application in geophysics,
including studies focused on fault
detection \cite{wu2019faultseg3d}, random noise attenuation \cite{saad2020deep}, and so on.
Some researchers have also attempted to design a DNN architecture that directly maps seismic
data to subsurface models in a fully data-driven manner.
Numerical experiments \cite{wu2019inversionnet,yang2019deep,zhu2023fourier,ding2022coupling,kazei2020deep} have demonstrated that the inverse operator for
full-waveform inversion can be acquired by training Convolutional Neural
Networks (CNNs) using datasets composed of wavefield and velocity model
pairs. Nonetheless,  it's important to note that the ill-posed and
complex nature of seismic inversion poses challenges for training CNNs
that are both sufficiently generalized and robust for inversion tasks.
This challenge is exacerbated by the need for extensive datasets containing
a diverse range of subsurface model/seismogram pairs.
In practical scenarios, obtaining such data can be a formidable task.

Compared to purely data-driven methods, FWI, as a substantial optimization problem,
exhibits enduring predictive capabilities and robust generalization.
The key to its prowess lies in its foundation on fundamental physical principles,
particularly those governing wave propagation. Therefore, it becomes imperative to
incorporate equation-based insights when harnessing deep learning techniques to
tackle the challenges posed by FWI. Rasht-Behesht et al. \cite{rasht2022physics} propose a new
approach to solve full-waveform inversion by Physics-Informed Neural Networks (PINNs),
it parameterizes the velocity model as a neural network and further jointly
optimizes the solutions of equations and the velocity model. Sun et al. \cite{sun2020theory}
develop a theory-guided system and choose an RNN as the most promising framework.
When the forward propagation of information through the network mirrors the forward
propagation of a wavefield through a heterogeneous medium, the associated RNN training
process can be viewed as a form of waveform inversion.  \cite{wu2019parametric, he2021reparameterized} parameterize the velocity parameters by a convolutional neural
network (CNN), such that the inversion amounts to reconstructing the weights of CNN. While these equation-informed unsupervised deep learning methods have made substantial advancements in the field of geophysical imaging, enhancing the numerical results of inversions, they still face challenges in accurately recovering the physical parameters of deeper subsurface layers, especially in complex geological structures. Therefore, the objective of this paper is to tackle the challenge and achieve accurate geophysical imaging by combining deep neural networks and FWI.

In this study, we introduce an unsupervised deep-learning approach designed
to tackle the complex problem of 2D full-waveform inversion.
This inverse problem is distinguished by its intrinsic nonlinearity and ill-posed characteristics,
making it a formidable challenge when employing traditional inversion methodologies.
Our approach is grounded in a re-parametrization technique for Bayesian inference,
leveraging a deep neural network with random weights.
By integrating likelihood functions contingent on physical equations with
prior velocity model information, we formulate the loss function used
for optimizing the neural network. To incorporate prior information into the neural network, we introduce a pretraining method that enables the neural network to learn an initial range of velocity parameter values, thus avoiding optimization from scratch during the unsupervised learning process.

Additionally, it is worth emphasizing that our method operates within an
unsupervised framework, eliminating the need for a dedicated training dataset.
Instead, it relies solely on an initial model and observed data.
To demonstrate the inversion capabilities of our algorithm,
we conduct an extensive series of numerical experiments.
These experiments provide compelling evidence that underscores
the feasibility and robustness of the proposed inversion methodology.
They convincingly demonstrate the ability of our approach to yield accurate
inversions for the inherently ill-posed inverse problem studied in this work,
effectively addressing the challenges posed by local minima during the full-waveform inversion process.

Main contributions are summarized as follows:
\begin{itemize}
	\item We propose an unsupervised deep learning approach that transforms the original grid-based reconstruction problem into an optimization problem for neural network parameters based on a physics equation and measurement data. Notably, the optimization of these network parameters requires no paired data whatsoever.
	\item By combining statistical inversion methods, neural network re-parameterization techniques, and variational Bayesian inference, we deduce an objective function for optimizing neural network parameters.
        \item We introduce a pretraining method to avoid initializing neural network parameters from scratch during the unsupervised learning process.
	\item We carry out an extensive series of numerical experiments, and the results from these experiments demonstrate that our approach surpasses the state-of-the-art inversion methods. In particular, our inversion method excels in accurately imaging complex subsurface structures and showcases robustness in the presence of measurement noise.
\end{itemize}

The paper is organized as follows. In Section \ref{sec2}, we introduce the problem setup of
full-waveform inversion.  Section \ref{sec3} describes the proposed method and algorithm, and
the neural network architecture and implementation details of our approach.
In Section \ref{sec4}, we experimentally compare our method against traditional FWI and a
deep learning method by different benchmark models.  Finally, we conclude our study in Section \ref{sec5}.

\section{Full waveform inversion}\label{sec2}
In this paper, we consider the acoustic wave equation, which characterizes
the propagation of pressure or sound waves in either fluid or solid media.
Specifically, we focus on the governing equation for acoustic waves in an
isotropic medium with uniform density, given by:
\begin{align}\label{2.1}
	\left \{
	\begin{aligned}
		&\frac{1}{v(x)^2}\frac{\partial^2u(x,t)}{\partial t^2}=\nabla^2u(x,t)+s(x,t;\xi),                                                                     \\
		& u(x,0)=0,                                                         \\
		& u_t(x,0)=0,
	\end{aligned}
	\right.
\end{align}
where $t$ is time, $x$ is the spatial location, $s$ is source function, and $v$ is the velocity
map for the subsurface medium. Here, we use a point source generated by a Ricker wavelet, represented as follows:
\begin{align*}
	s(x,t;\xi)=s_0(t)\delta(x-\xi),
\end{align*}
where $s_0(t)=(1-2\pi^2 f^2t^2)e^{-\pi^2f^2t^2}$
is the amplitude of the Ricker wavelet with frequency $f$, $\delta$ denotes the Dirac delta
function, and the $\xi$ is the predetermined horizontal location of the source.
To simulate a realistic application, we applied the reflection boundary condition to the top
surface of the physical domain. For the remaining boundaries, we implement absorbing boundary
layers to effectively simulate wave propagation in an unbounded medium.
To achieve this, we employ a perfectly matched layer (PML) approach \cite{Komatitsch+Tromp+2003},
which provides absorbing boundary conditions along the left, right, and bottom sides.
To show the dependence of the solution $u$ on the parameters $v,~\xi$, we will use $u(x,t;v,\xi)$ to represent the solution to Equation (\ref{2.1}).

\begin{figure}[H]
	\centering
	\subfigure{
		\begin{minipage}[t]{1\linewidth}
			\centering
			\includegraphics[width=3.5in]{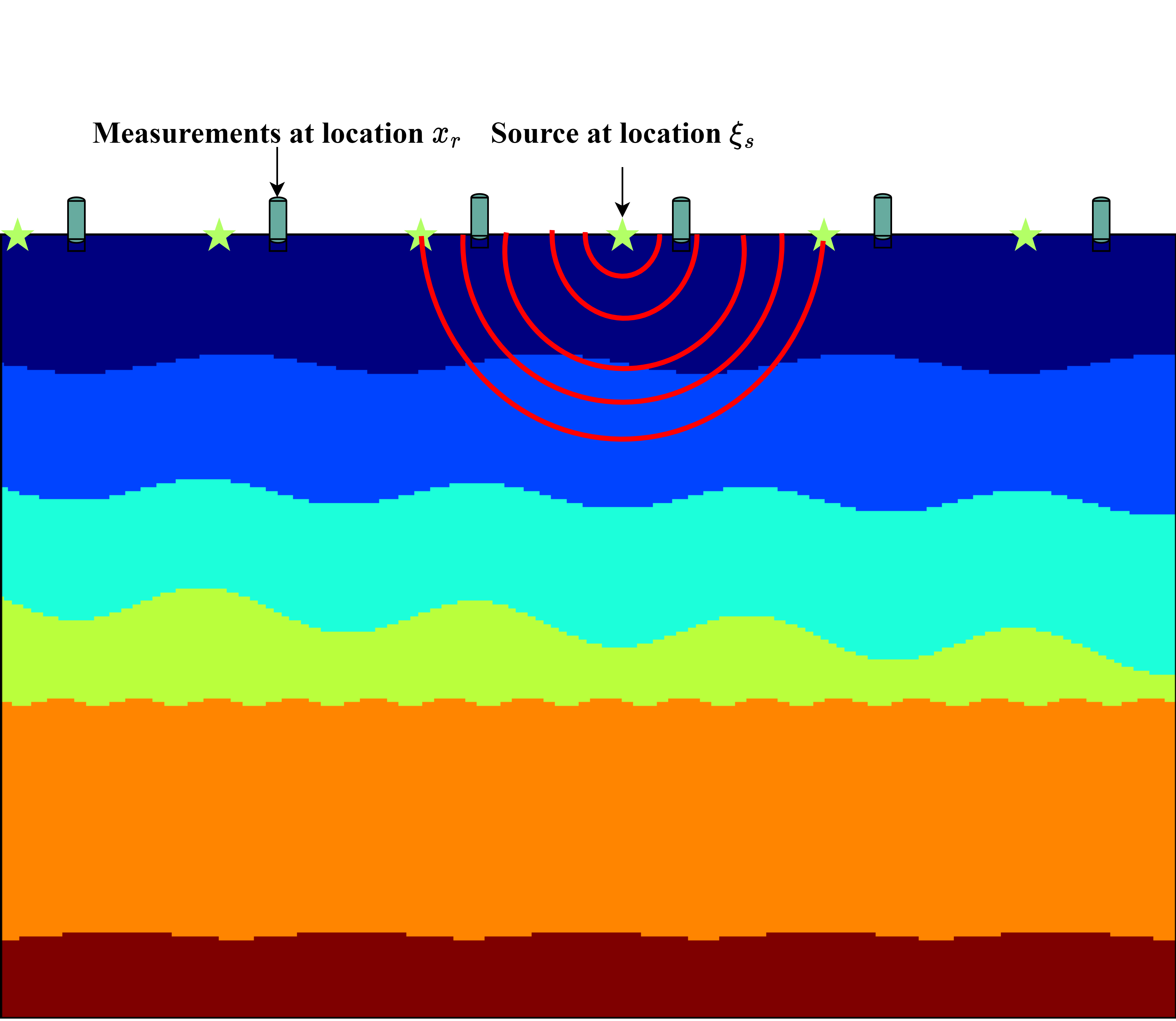}
		\end{minipage}
	}%
	\centering
	\caption{A conceptual illustration of full-waveform inversion.  The signal emitted by the sources (Star) is recorded at the
sensors (Cylinder).}\label{fig_wave_source}
\end{figure}

Full waveform inversion is primarily concerned with the recovery of the subsurface velocity,
denoted as $v$, from the observed wavefield, represented by $u$. In practical applications,
we face limitations on both the number of available receivers and their placement, typically constrained to the surface.
To achieve an accurate reconstruction of the subsurface velocity, it is common practice
to employ multiple sources to generate the surface wavefield. Specifically,
we consider a scenario with $n_s$ sources, denoted as $s(x, t; \xi_s)$,
where $s=1,\cdots,n_s$.
For each source function $s(x,t;\xi_s)$, we denote the corresponding wavefield as $u(x,t;v,\xi_s)$.
The observed wavefield is collected only at the surface, with a total of $n_r$ receivers positioned
at locations denoted as $x_r$, where $r=1,\cdots,n_r$. A conceptual
illustration of full-waveform inversion is provided in Figure \ref{fig_wave_source}.

In this paper, we define our computational domain as $\Omega = [0, L_x] \times [0, L_z] \subset \mathbb{R}^2$.
We discretize this domain and represent the discrete coordinates as $(x_i,z_j)~(i=1,\cdots,m_x,~j=1,\cdots,m_z)$.
We denote $v=\{v(x_i, z_j)\}_{i=1,\cdots,m_x,~j=1,\cdots,m_z}$, i.e., $v\in \mathbb{R}^{m_x \times m_z}$.
Given a velocity model $v$ and a source function $s(x, t; \xi_s)$,
we solve equation (\ref{2.1}) to compute the observed wavefield $u_s(x_r, t; v, \xi_s) = Ru(x, t; v, \xi_s)$,
where the operator $R$ calculates the wavefield $u(x, t; v, \xi_s)$ at a receiver location $x_r$.
Consequently, we obtain the observed data $d_s(x_r, t) = u_s(x_r, t; v, \xi_s) + \eta_s$,
with $\eta_s\sim N(0, \Sigma)$, where $\Sigma = \sigma^2 I$.
Thus, the goal of $l_2$ norm full-waveform inversion is to recover the velocity $v$ by minimizing
the following cost function:
\begin{align*}
	\min\limits_{v} L(v)=\min\limits_{v}\{\sum_{s=1}^{n_s}\sum_{r=1}^{n_r}\|u_s(x_r,t;v,\xi_s)-d_s(x_r,t)\|^2\}.
\end{align*}
Additionally, we introduce the forward operator $F(v)$, defined as:
\begin{align*}
	F(v)=u,
\end{align*}
where $v\in \mathbb{R}^{m_x \times m_z}$, $u\in \mathbb{R}^{n_s\times n_r\times n_T}$.
We can then express the general formulation of the observed data as:
\begin{align}\label{2.2}
	d=F(v)+\eta,
\end{align}
where $d,~\eta\in\mathbb{R}^{n_s\times n_r\times n_T}$, $\eta\sim N(0,\Sigma)$ with $\Sigma=\sigma^2 I$.

In general, full-waveform inversion based on the $l_2$-norm relies heavily on the accuracy of the initial velocity model, primarily due to the profound nonlinearity of FWI. The optimization process of $l_2$-norm FWI is susceptible to getting stuck in local minima caused by the cycle skipping problem.
To address this issue, numerous researchers have explored alternative loss functions, such as $l_1$ norm, Huber and hybrid $l_1-l_2$ norm \cite{brossier2010data}. Additionally, the use of the Wasserstein metric, based on optimal transport, has shown promise in FWI and yielded competitive results. Nevertheless, despite these advancements, obtaining accurate inversion results for complex subsurface structures remains an ongoing challenge. Therefore, there is a compelling need to investigate more accurate numerical inversion methods aimed at addressing the challenge of solving the full-waveform inversion problem in complex underground structures.

\section{Our approach}\label{sec3}
In this section, we will provide detailed explanations of our approach for FWI with noisy measurements. First, by leveraging the Bayesian inversion framework and neural network re-parameterization techniques, we introduce the rationale and motivation behind this study.  Secondly, by variational Bayesian inference, we formally deduce a loss function that used to optimize the neural network parameters. Thirdly, to encode the prior information of the velocity parameters into the neural network, we propose a pretraining method for the neural network. Subsequently, using the trained neural network, we propose a conditional
mean estimator for estimating velocity parameters $v$.
Finally, we provide details regarding the structure of the neural network used to parameterize the velocity parameter $v$, along with experimental specifics.
\subsection{Motivations}
In this paper, we propose an unsupervised deep learning approach to tackle the intricate nonlinear inverse problem.
Our method is grounded in a re-parametrization technique using Deep Neural Networks (DNNs)
for Bayesian inference.
Within the Bayesian framework, we encapsulate our prior beliefs about the unknown velocity,
denoted as $v$ in the prior probability distribution, which is represented as $p(v)$.
In addition, we assume that the unknown parameter $v$
is statistically independent of the noise $\eta$.
Then, the distribution of the $v$ conditioned on the measurement data $d$, i.e.,
the posterior distribution $p(v|d)$ follows the Bayes' rule:
\begin{align}\label{3.1}
	p(v|d)=\frac{1}{p(d)}p(d|v)p(v),
\end{align}
where $p(d|v)$ denotes the likelihood
\begin{align*}
	p(d|v)\propto \exp(-\|F(v)-d\|_{\Sigma}^2),
\end{align*}
with
\begin{align*}
	\|F(v)-d\|_{\Sigma}^2=(F(v)-d)^{T}\Sigma^{-1}(F(v)-d).
\end{align*}

In conventional Bayesian inversion methods, the need for computational
feasibility often leads to modeling the prior distribution $p(v)$ by a simple distribution. The oversimplified
prior distribution has prompted us to
explore alternative ways of representing either the prior distribution $p(v)$
or the posterior distribution $p(v|d)$.
Drawing inspiration from the potent representational capabilities of deep neural networks (DNNs),
we propose employing a re-parametrization technique for the Bayesian inverse problem.
This technique involves re-expressing the variable $v$ using a DNN with random weights:
\begin{align}\label{3.1.1}
	v=m(z_0;\bm{\theta}),
\end{align}
where $z_0$ represents a fixed random tensor with the size $[1,1,m_z,m_x]$, and $\bm{\theta}$ comprises random variable.
It is noteworthy that, through re-parametrization, the prior distribution $p(v)$
can become highly intricate, even when the random variable $\bm{\theta}$ is
assumed by simple distributions.

Following the re-parametrization, the variable for Bayesian inference become
the random weight denoted as $\bm{\theta}$ associated with the neural network.
Instead of inferring the value of $v$ from the posterior distribution (\ref{3.1}),
our objective now shifts to inferring the variable $\bm{\theta}$ from the posterior
distribution $p(\bm{\theta}|d)$, which is defined as:
\begin{align}\label{3.2}
	p(\bm{\theta}|d)\propto p(d|\bm{\theta})p(\bm{\theta}),
\end{align}
where
\begin{align}\label{3.2.1}
	p(d|\bm{\theta})\propto \exp(-\|F(m(z_0;\bm{\theta}))-d\|_{\Sigma}^2),
\end{align}
and $p(\bm{\theta})$ represents the prior distribution for the random weight $\bm{\theta}$.

Due to the infeasibility of high-dimensional computations in solving the posterior distribution $p(\bm{\theta}|d)$ as outlined in (\ref{3.2}), approximate inference methods have been developed. The central challenge lies in creating an expressive approximation to the true posterior while maintaining computational efficiency and scalability, especially within modern deep-learning architectures. Variational inference stands out as a popular deterministic approximation approach for addressing this challenge. In variational inference, the approximate posterior is assumed to be fully factorized distributions, often referred to as mean-field variational inference \cite{graves2011practical,kingma2016improved}. In general, mean-field variational approximation promotes computational tractability and effective optimization. However, it does have limitations in capturing the intricate structure of the true posterior \cite{nguyen2021structured}.

In recent years, there have been numerous studies demonstrating the use of Dropout regularization as an approximation in Bayesian inference models \cite{gal2016dropout,kingma2015variational,nguyen2021structured}. The literature highlights that inference methods based on Bayesian Dropout have delivered competitive performance in predictive accuracy across various tasks, offering several promising avenues for enhancing approximate inference in Bayesian models.

In this paper, we employ the variational Dropout approximation method \cite{gal2016dropout}, which approximates the posterior distribution $p(\bm{\theta}|d)$ using a set of approximation distributions $q(\bm{\theta}|\bm{\mu})$ parametrized by $\bm{\mu}$, i.e., we introduce a set of distributions characterized by the following relationship:
\begin{equation}\label{3.2.2}
	\bm{\theta}=\bm{\mu}\odot \bm{b}:~\theta_{i}=\mu_{i}*b_{i},~i=1,\cdots,N,
\end{equation}
where $\mu_i$ represents the distribution parameter associated
with $\theta_i$, and $b_i\sim \bm{B}(p_i)$ follows a Bernoulli
distribution with a probability $p_i$.
In other words, the probability density function of $b_i$ is defined as:
\begin{align*}
	p(b_i)=p_i^{b_i}(1-p_i)^{1-b_i},~b_i\in \{0,1\}.
\end{align*}

\begin{figure}[H]
	\centering
	\subfigure{
		\begin{minipage}[t]{1\linewidth}
			\centering
			\includegraphics[width=6.5in]{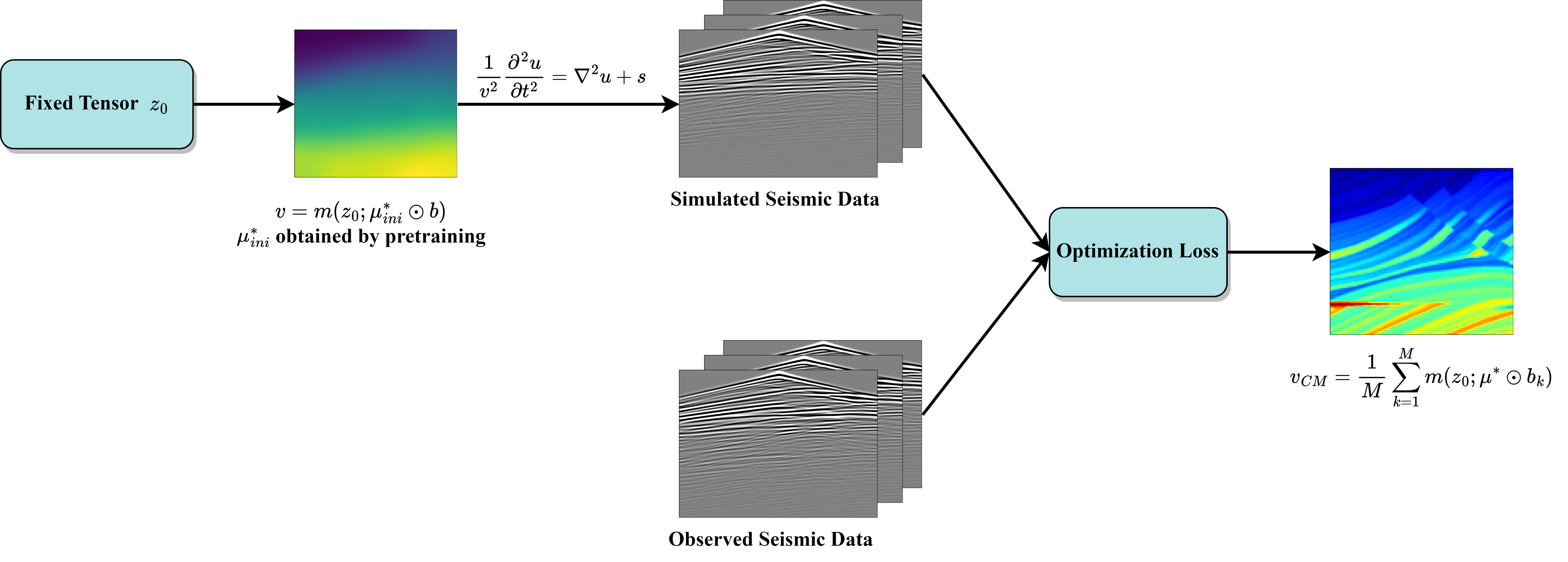}
		\end{minipage}
	}%
	\centering
	\caption{Schematic workflow of our proposed method for full-waveform inverion. $\bm{\mu}_{ini}^{*}$ is obtained through the pretraining detailed in the following sections.}\label{fig_digram}
\end{figure}

 The optimal parameters $\bm{\mu}$ of the distribution $q(\bm{\theta}|\bm{\mu})$ are
chosen based on a metric that quantifies the difference between the approximation
and the actual posterior distribution. A commonly used metric in variational
inference is the Kullback-Leibler (KL) divergence.
In general, the KL divergence can be expressed as:
\begin{align*}
	D_{KL}(q(x)\|p(x)) = \int q(x)\log\frac{q(x)}{p(x)}dx.
\end{align*}
Thus, in this paper, we choose the optimal parameters $\bm{\mu}$ by solving the following
minimization problem:
\begin{align*}
		\mathop{\min}_{\bm{\mu}}D_{KL}(q(\bm{\theta}|\bm{\mu})\|p(\bm{\theta}|d)).
\end{align*}
Through some detailed calculations (Specific computation details will be deferred to the subsequent context.), we can deduce the following optimization problem:
\begin{align}\label{3.2.3}
	\mathop{\min}_{\bm{\mu}}L_{W_1}(\bm{\mu})&=\mathop{\min}_{\bm{\mu}}E_{\bm{b}\sim \bm{B}(p)}\big\{\sum_{s=1}^{n_s}\sum_{r=1}^{n_r}\|u_s(x_r,t;m(z_0;\bm{\mu}\odot\bm{b}),\xi_s)-d_s(x_r,t)
	\|_{W_1}+\alpha TV(m(z_0;\bm{\mu}\odot\bm{b}))\big\},
\end{align}
where $\|\cdot\|_{W_1}$ denotes $W_1$ distance \cite{zhang2022optimal}, $TV$ represents the Total Variation regularization and $\alpha$ is a regularization
parameter. Then, we denote the reconstructed velocity model by
\begin{align}\label{3.2.4}
	v_{CM} \approx \frac{1}{M}\sum_{k=1}^{M}m(z_0;\bm{\mu}^*\odot \bm{b_k}),~ \bm{b}_k\sim \bm{B}(p),
\end{align}
where $\bm{\mu}^*$ is a minimizer of problem (\ref{3.2.3}). In summary, a schematic diagram of our approach is presented in Figure \ref{fig_digram}.

Next, we will provide a detailed explanation of how to obtain the optimization problem (\ref{3.2.3}) and how to derive the estimation formula for the reconstructed velocity model (\ref{3.2.4}).

\subsection{Constructing the loss}

As we employ $q(\bm{\theta}|\bm{\mu})$ to approximate the posterior
distribution $p(\bm{\theta}|d)$, we determine the optimal approximation
by minimizing the following KL-divergence:
\begin{align}\label{3.3}\nonumber
	\mathop{\min}_{\bm{\mu}}D_{KL}(q(\bm{\theta}|\bm{\mu})\|p(\bm{\theta}|d))
	&=\mathop{\min}_{\bm{\mu}}\int q(\bm{\theta}|\bm{\mu})\log\frac{q(\bm{\theta}|\bm{\mu})}{p(\bm{\theta}|d)}d\bm{\theta}\\ \nonumber
	%&=\mathop{\min}_{\bm{\mu}}\int q(\bm{\theta}|\bm{\mu})\log\frac{q(\bm{\theta}|\bm{\mu})}{\frac{1}{p(d)}p(d|\bm{\theta})p(\bm{\theta)}}d\bm{\theta}\\ \nonumber
	%&=\mathop{\min}_{\bm{\mu}}\int q(\bm{\theta}|\bm{\mu})\log\frac{q(\bm{\theta}|\bm{\mu})}{p(d|\bm{\theta})p(\bm{\theta})}d\bm{\theta}
	%+\int q(\bm{\theta}|\bm{\mu})\log p(d)d\bm{\theta}\\ \nonumber
	&\propto\mathop{\min}_{\bm{\mu}}[\int q(\bm{\theta}|\bm{\mu})\log\frac{q(\bm{\theta}|\bm{\mu})}{p(\bm{\theta})}d\bm{\theta}
	-\int q(\bm{\theta}|\bm{\mu})\log p(d|\bm{\theta})d\bm{\theta}]\\
	&=\mathop{\min}_{\bm{\mu}}[D_{KL}(q(\bm{\theta}|\bm{\mu})\|p(\bm{\theta}))-E_{\bm{\theta}\sim q(\bm{\theta}|\bm{\mu})}\log p(d|\bm{\theta})].
\end{align}

For the first term in formula (\ref{3.3}), assuming that $p(\bm{\theta})$ follows
a uniform distribution over a sufficiently large region,
we can represent $p(\theta_i)$ as $1/s_i$, where $s_i$ corresponds
to the length of the domain associated with the variable $\theta_i$.
By (\ref{3.2.2}), we have $q(\theta_i|\mu_i)=p_i^{\frac{\theta_i}{\mu_i}}(1-p_i)^{1-\frac{\theta_i}{\mu_i}}$, $\theta_i\in\{0,\mu_i\}$, where we abuse
the notion $\frac{0}{0}=1$.
Consequently, we obtain:
\begin{align*}
	D_{KL}(q(\bm{\theta}|\bm{\mu})\|p(\bm{\theta}))&=\sum_{i=1}^{N}D_{KL}(q(\theta_i|\mu_i)\|p(\theta_i))\\
	&=\sum_{i=1}^{N}\int q(\theta_i|\mu_i)\log \frac{q(\theta_i|\mu_i)}{p(\theta_i)}d\theta_i\\
	&=\sum_{i=1}^{N}\big(\int q(\theta_i|\mu_i)\log q(\theta_i|\mu_i)d\theta_i+\int q(\theta_i|\mu_i)\log s_id\theta_i\big)\\
	&=\sum_{i=1}^{N}\big((1-p_i)\log(1-p_i)+p_i\log p_i+\log s_i\big)=C_0,
\end{align*}
where the constant $C_0$ independent of the parameters $\bm{\mu}$.

By (\ref{3.2.1}), we have
\begin{align*}
	\log p(d|\bm{\theta})\propto -\|F(m(z_0;\bm{\theta}))-d\|_{\Sigma}^2,
\end{align*}
then
\begin{align*}
	\mathop{\min}_{\bm{\mu}}D_{KL}(q(\bm{\theta}|\bm{\mu})\|p(\bm{\theta}|d))
	\propto\mathop{\min}_{\bm{\mu}}E_{\bm{\theta}\sim q(\bm{\theta}|\bm{\mu})}\|F(m(z_0;\bm{\theta}))-d\|_{\Sigma}^2.
\end{align*}
Using (\ref{3.2.2}), we can deduce that
\begin{align*}
	\mathop{\min}_{\mu}E_{\bm{\theta}\sim q(\bm{\theta}|\bm{\mu})}\|F(m(z_0;\bm{\theta}))-d\|_{\Sigma}^2
	&=\mathop{\min}_{\bm{\mu}}\int\|F(m(z_0;\bm{\theta}))-d\|_{\Sigma}^2 q(\bm{\theta}|\bm{\mu})d\bm{\theta}\\
	&=\mathop{\min}_{\bm{\mu}}\int \|F(m(z_0;\bm{\mu}\odot \bm{b}))-d\|_{\Sigma}^2 \bm{B}(p)d\bm{b}\\
	&=\mathop{\min}_{\bm{\mu}}E_{\bm{b}\sim \bm{B}(p)}\|F(m(z_0;\bm{\mu}\odot\bm{b}))-d\|_{\Sigma}^2.
\end{align*}

To mitigate the risk of potential overfitting,
we introduce an additional regularization term to the estimation process.
We adopt the commonly employed Total Variation (TV) regularization and
incorporate it into the loss function. Therefore,
the loss function used for training the neural network becomes:
\begin{align*}
	\mathop{\min}_{\bm{\mu}}L(\bm{\mu})&=\mathop{\min}_{\bm{\mu}}E_{\bm{b}\sim \bm{B}(p)}\big\{\|F(m(z_0;\bm{\mu}\odot\bm{b}))-d\|_{\Sigma}^2+\alpha TV(m(z_0;\bm{\mu}\odot\bm{b}))\big\}\\
	&=\mathop{\min}_{\bm{\mu}}E_{\bm{b}\sim \bm{B}(p)}\big\{\sum_{s=1}^{n_s}\sum_{r=1}^{n_r}\|u_s(x_r,t;m(z_0;\bm{\mu}\odot\bm{b}),\xi_s)-d_s(x_r,t)
	\|^2+\alpha TV(m(z_0;\bm{\mu}\odot\bm{b}))\big\}.
\end{align*}

In general, full-waveform inversion using the least-squares
loss can encounter challenges like getting trapped in local
minima due to issues such as cycle skipping in wave oscillations.
In contrast, the Wasserstein distance, rooted in optimal transport theory \cite{villani2009optimal},
is convex with respect to shifted patterns when
used to positive functions, which is an advantageous
property for solving FWI. Here, we consider $W_1$ distance \cite{zhang2022optimal} as the measure to quantify
the discrepancy between the observed and simulated data. Consequently, the corresponding loss function is formulated as:
\begin{align}\label{3.4}
	\mathop{\min}_{\bm{\mu}}L_{W_1}(\bm{\mu})&=\mathop{\min}_{\bm{\mu}}E_{\bm{b}\sim \bm{B}(p)}\big\{\sum_{s=1}^{n_s}\sum_{r=1}^{n_r}\|u_s(x_r,t;m(z_0;\bm{\mu}\odot\bm{b}),\xi_s)-d_s(x_r,t)
	\|_{W_1}+\alpha TV(m(z_0;\bm{\mu}\odot\bm{b}))\big\}.
\end{align}

\subsection{Encoding the prior}
In the field of inverse problem research, the accuracy of the reconstructed solution is heavily reliant on the amount of prior information embedded during the inversion process. This prior information may encompass the smoothness, sparsity, or the choice of initial values for the inverse problem solution. Generally, the more known prior information is available, the more accurate the resulting inverse problem solution will be. Conversely, when there is limited prior knowledge about the inversion parameters, obtaining a high-precision solution for the inverse problem becomes more challenging.

In our research, full-waveform inversion is a well-known nonlinear inverse problem, which leads to a minimization problem with multiple local minima. Therefore, in order to achieve a more precise solution to this inverse problem, it is necessary for the neural network to learn some prior information.
In our implementation, we empower the DNN, denoted as $m(z_0;\bm{\theta})=m(z_0;\bm{\mu}\odot\bm{b})$,
to learn some prior information regarding the initial value $v_{ini}$.
To accomplish this, we minimize the following problem:
\begin{align}\label{3.5}
	\bm{\mu}_{ini}^{*} = \arg\min J(\bm{\mu}) = \arg\min E_{\bm{b}\sim\bm{B}(p)}\|m(z_0;\bm{\mu}\odot\bm{b})-v_{ini}\|_{l_1},
\end{align}
where, $\|\cdot\|_{l_1}$ denotes the $l_1$ norm, and $v_{ini}$ represents an initial model of the velocity parameter $v$.

\subsection{Inference}
After training the neural network by minimizing the loss function as specified in (\ref{3.4}),
we can obtain an approximation for the posterior distribution $p(\bm{\theta}|d)$, which is represented as $q(\bm{\theta}|\bm{\mu}^*)$.
During the testing phase of our method, we estimate the velocity parameter $v$
using a conditional mean estimator. Given the observation data $d$,
the velocity estimator $v_{CM}$ is defined as:
\begin{align*}
	v_{CM} = \int vp(v|d)dv.
\end{align*}
Utilizing (\ref{3.1.1}), we can deduce that:
\begin{align*}
	v_{CM} = \int m(z_0;\bm{\theta})p(\bm{\theta}|d)d\bm{\theta}.
\end{align*}
By employing the variational inference method,
where the distribution $q(\bm{\theta}|\bm{\mu}^{*})$ approximates the posterior distribution $p(\bm{\theta}|d)$, then we can express:
\begin{align*}
	v_{CM} \approx \int m(z_0;\bm{\theta})q(\bm{\theta}|\bm{\mu}^{*})d\bm{\theta} = \int m(z_0;\bm{\mu}^*\odot\bm{b})\bm{B}(p)d\bm{b}.
\end{align*}
In practice, the integration is computed using Monte Carlo (MC) method.
After the neural network is trained, we estimate the velocity parameter $v$ by
\begin{align}\label{3.6}
	v_{CM} \approx \frac{1}{M}\sum_{k=1}^{M}m(z_0;\bm{\theta}_k) = \frac{1}{M}\sum_{k=1}^{M}m(z_0;\bm{\mu}^*\odot \bm{b_k}),
\end{align}
where $\bm{b}_k\sim \bm{B}(p)$.

Based on the description provided above for our inversion method, we summarize the algorithm in Algorithm 1.
\begin{algorithm}[H]
	\caption{Our algorithm}
	\begin{algorithmic}[1]
		\STATE \textbf{Inputs:} $z_0$, $v_{ini}$, $\epsilon$, $It_{pmax}$, $It_{Imax}$;\\
		\STATE \textbf{Phase 1:}\\
		Solve the following problem\\
		\begin{align*}
			\min J(\bm{\mu}),
		\end{align*}
		where $J(\bm{\mu})$ defined by (\ref{3.5}).
		The iterations stops and outputs $\bm{\mu}_{ini}^{*}$ when $J(\bm{\mu}_{ini}^{*})<\epsilon$ or reaches the maximum number of iterations $It_{pmax}$.
		
		\STATE \textbf{Phase 2:}\\
		Let $\bm{\mu}_0=\bm{\mu}_{ini}^{*}$ and solve the minimization problem (\ref{3.4}), where $\bm{\mu}_0$ is a initial value of parameter $\bm{\mu}$ of the minimization problem (\ref{3.4}). The iterations stop when reach the maximum
		number of iterations $It_{Imax}$.
		
		\STATE \textbf{Inference:}\\
		Calculate $v_{CM}$ by (\ref{3.6}).
	\end{algorithmic}
\end{algorithm}

\begin{comment}
    In this section, we have provided a comprehensive exposition of our proposed inversion method. In fact, we approximate the posterior distribution $p(\bm{\theta}|d)$ using a neural network with dropout $q(\bm{\theta}|\bm{\mu})$. Dropout, along with its variations in neural networks, can be viewed as a Bayesian approximation \cite{gal2016dropout} of a well-established probability model: the Gaussian process. Once the training of the neural network with dropout is completed, we employ a Monte Carlo-based approximation of the conditional mean estimator of velocity parameters to represent the reconstructed solution of the full waveform inversion.
\end{comment}
\subsection{Implementation details}

\begin{figure}[H]
	\centering
	\subfigure[The overall architecture of the proposed encoder-decoder neural network.]{
		\begin{minipage}[t]{1\linewidth}
			\centering
			\includegraphics[width=3.5in]{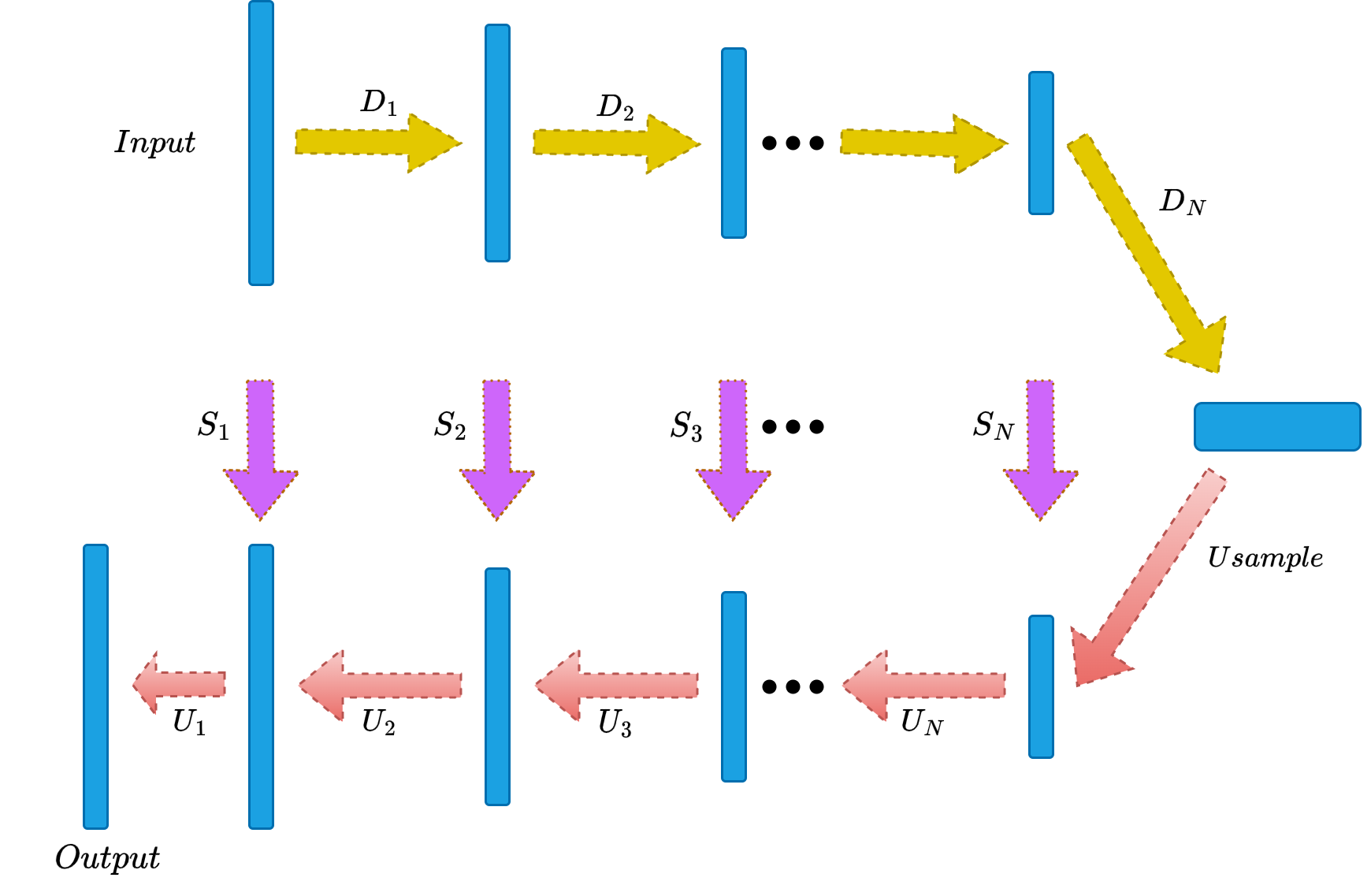}
		\end{minipage}
	}%
	
	\subfigure[Each component block of the proposed encoder-decoder neural network.]{
		\begin{minipage}[t]{1\linewidth}
			\centering
			\includegraphics[width=3.5in]{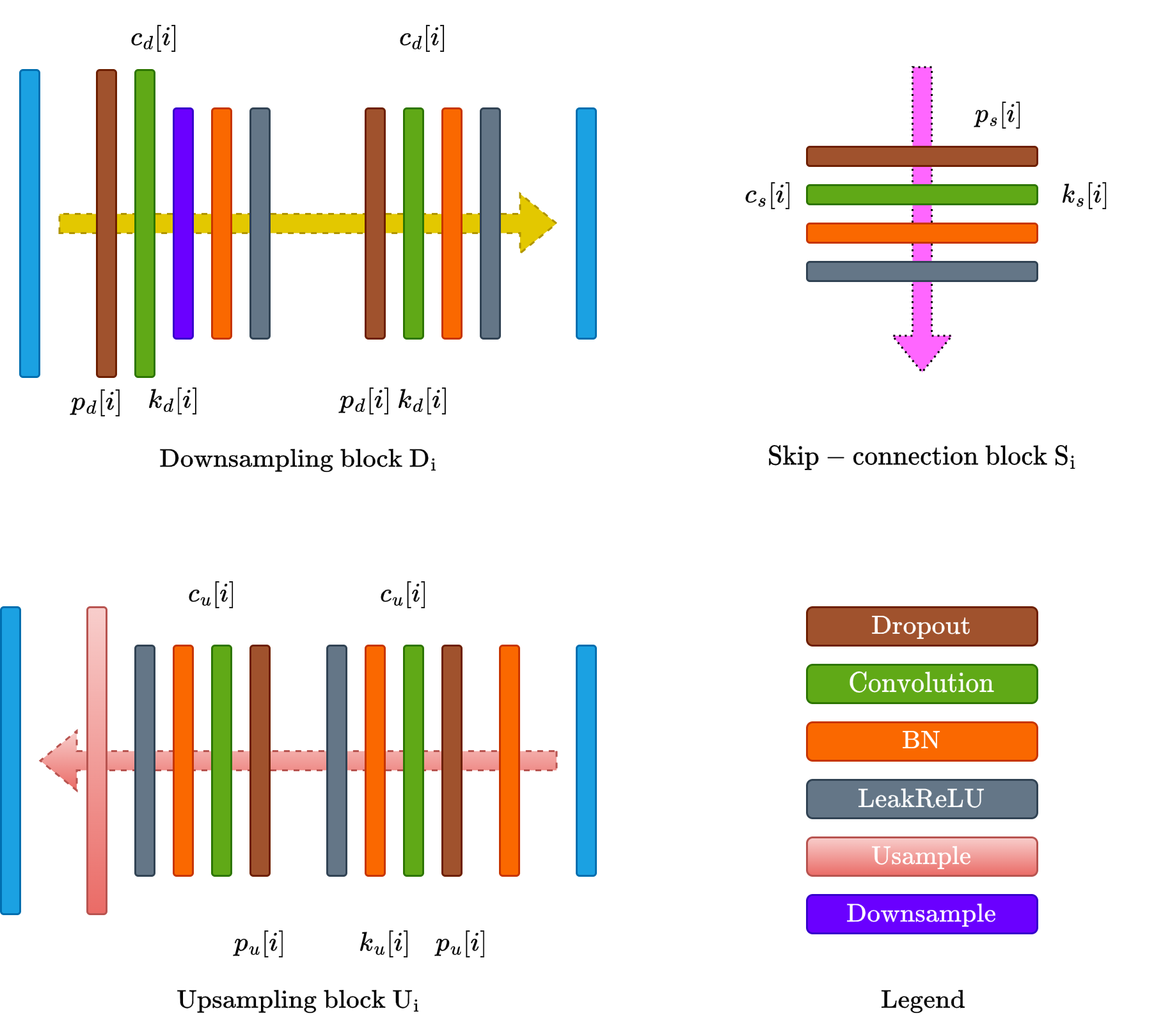}
		\end{minipage}
	}%
	\centering
	\caption{Diagram of the network used for evaluating the proposed method.}\label{fig_1}
\end{figure}

In this paper,
we employ an encoder-decoder neural network \cite{ulyanov2018deep} that incorporates skip-connections to parameterize the velocity model $v$.
Figure \ref{fig_1} (a) provides a detailed illustration of the network architecture.
Within this network diagram, we represent the components of the encoder-decoder
architecture as follows: downsampling by $D_i$, upsampling by $U_i$, and
skip-connection by $S_i$. Each block of the encoder-decoder has specific
parameters denoted as follows: the number of filters at depth $i$ for
downsampling, upsampling, and skip-connections as $c_d[i]$, $c_u[i]$, and $c_s[i]$, respectively.
The kernel sizes are represented by $k_d[i]$, $k_u[i]$, $k_s[i]$, and the dropout
probabilities for downsampling, upsampling, and skip-connections are $p_d[i]$, $p_u[i]$, and $p_s[i]$, respectively.

In our experiments, we set $N=5$, $c_d[i]=c_u[i]=128$, $c_s[i]=4$,
and use $k_d[i]=k_u[i]=3$ and $k_s[i]=1$ for the convolution layers.
We employ the Leaky ReLU as the non-linear activation function with
a slope of 0.1. Downsampling layers are implemented through convolutional
blocks with a stride of 2, while bi-linear interpolation is used for upsampling.
We set the dropout probability to $p_d[i]=p_u[i]=0$ for both $D_i$ and $U_i$ and to $p_{s}[i]=0.3$ for $S_i$.

\section{Numerical experiments}\label{sec4}
In this section, we demonstrate the advantages of our method through three numerical tests,
which include the Marmousi model, the Overthrust model, and the Marmousi2 model.
To evaluate the performance of our proposed method, we will use the following
relative $l_2$ error, defined as:
\begin{align*}
	\text{Relative}~l_2~\text{error}=\frac{\|\tilde{v}-v\|}{\|v\|},
\end{align*}
where $\tilde{v},~v$ represent the inverted velocity model and the ground truth model, respectively.
Additionally, we calculate the Signal-to-Noise Ratio (SNR) between the inverted velocity model $\tilde{v}$ and the ground truth model $v$ using the formula:
\begin{align*}
	\text{SNR}(\tilde{v}, v)=10\log_{10}\frac{\|v\|^2}{\|\tilde{v}-v\|^2}.
\end{align*}
Furthermore, to further assess the similarity between the
inverted velocity $\tilde{v}$ and the ground truth $v$,
we compute the Structural Similarity Index (SSIM) defined as:
\begin{align*}
	\text{SSIM}(\tilde{v},v)=\frac{(2\mu_{v}\mu_{\tilde{v}}+c_1)(2\sigma_{v\tilde{v}}+c_2)}{(\mu_{v}^2+\mu_{\tilde{v}}^2+c_1)(\sigma_v^2+\sigma_{\tilde{v}}^2+c_2)},
\end{align*}
where $\mu_v$, $\mu_{\tilde{v}}$, $\sigma_v$, $\sigma_{\tilde{v}}$, and $\sigma_{v\tilde{v}}$ are the local mean, standard deviation, and cross-covariances for $v$ and $\tilde{v}$,
respectively. Higher SSIM and SNR values indicate a better quality of the inversion result.

In this study, to demonstrate the effectiveness of our proposed method,
we conduct a comparison with state-of-the-art FWI methods.
This includes two traditional inversion methods utilizing $l_2$ loss (FWI($l_2$))
and the $1D-W_1$ loss (FWI($W_1$)) \cite{zhang2022optimal}, as well as a deep learning-based method (DNN-FWI) \cite{he2021reparameterized}.

To simulate wave propagation for Equation (\ref{2.1}),
we use a uniform grid and a finite difference scheme with
second-order accuracy in the time domain and eighth-order
accuracy in the spatial domain.
We choose a Ricker wavelet with a peak frequency of 5 Hz as the source wavelet.
The time interval for forward simulation is set to 3 ms,
and the total recording time is configured as 6 s.
We position thirty equally spaced sources on the surface at $z=0$
and take the number of receivers to match the model size $m_x$ for different velocity models.

To facilitate these simulations, we leverage the Deepwave toolbox \cite{richardson_alan_2023},
which offers wave propagation modules and facilitates automatic updates of variables.
To ensure a fair comparison, we fine-tune the hyper-parameters and implement each
method multiple times to obtain the best solutions.

Our method is implemented in the PyTorch interface on an NVIDIA 3090 GPU  graphics card with 24G memory.
In Algorithm 1, we set $It_{pmax}=5000$, $It_{Imax}=1000$, and $\epsilon=1\times 10^{-3}$.
In phase 1 of Algorithm 1, we apply the Adam optimizer with a learning rate of 0.01 to solve the minimization problem (\ref{3.5}).
In phase 2 of Algorithm 1, we minimize the cost function (\ref{3.4}) by Adam
optimizer with learning rate  $5\times 10^{-4}$.
We initialize the value of the fixed random tensor $z_0$ in (\ref{3.1.1}) as a uniform distribution
with a size of $[1,1,m_z,m_x]$.
 We take the regularization parameter $\alpha=6\times 10^{-7}$ in (\ref{3.4}) and The number of Monte Carlo samples $M=50$ in (\ref{3.6}).
 To ensure a fair comparison,
we set the number of iterations for all methods to be 1000.

\subsection{Marmousi model}
\begin{figure}[htbp]
	\centering
	\includegraphics[width=1.0\textwidth]{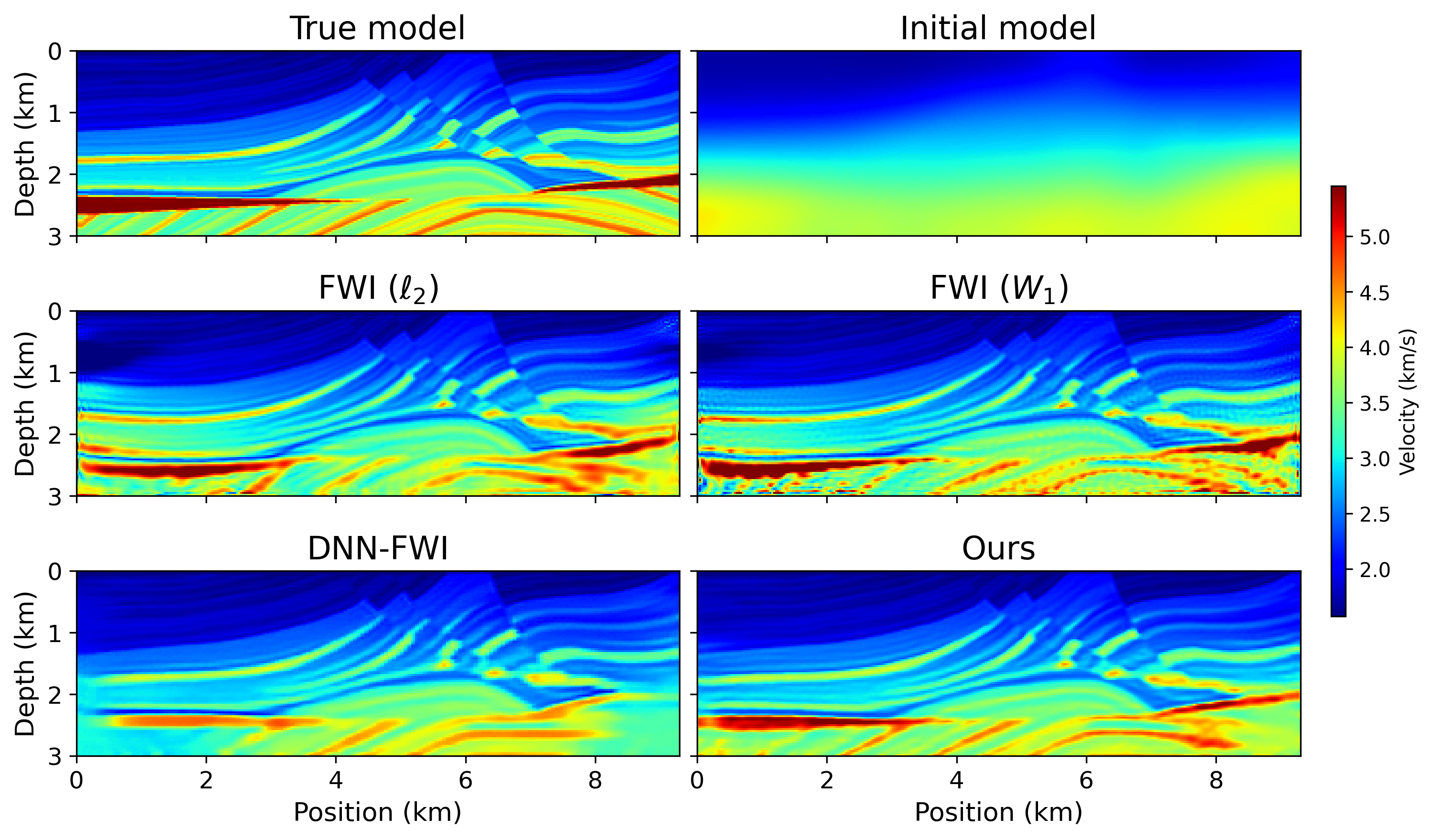}
	\caption{Inversion results for the noise-free observed data of the Marmousi model.
		The top row displays the true velocity model and the initial guess.
		The second row shows the inverted results using traditional FWI constrained by the $l_2$-norm and the $1D-W_1$ distance.
		The third row presents the inverted results obtained through DNN-FWI and our method.}
	\label{fig_mar0}
\end{figure}

In this part, we validate the effectiveness of our proposed method using
the Marmousi model \cite{versteeg1994marmousi}. The Marmousi model is a well-known
acoustic velocity model and serves as a standard benchmark model in geophysical exploration.
For our experiments, we reduce the resolution of the velocity model
to $(z\times x) = (100 \times 310)$ with a spatial grid increment of $0.03$ km.
The velocity values in the model range from $1472$ m/s to $5772$ m/s and the true model is depicted in Figure \ref{fig_mar0}.

Figure \ref{fig_mar0} presents a comparison of subsurface velocity models reconstructed
through various inversion methods using noise-free measurements. These methods yield reasonable
results overall. However, FWI($l_2$), FWI($W_1$), and DNN-FWI demonstrate limited accuracy
in deeper regions, contrasting with our approach, which excels in these areas.

\begin{figure}[H]
	\centering
	\subfigure{
		\begin{minipage}[t]{0.5\linewidth}
			\centering
			\includegraphics[width=3in]{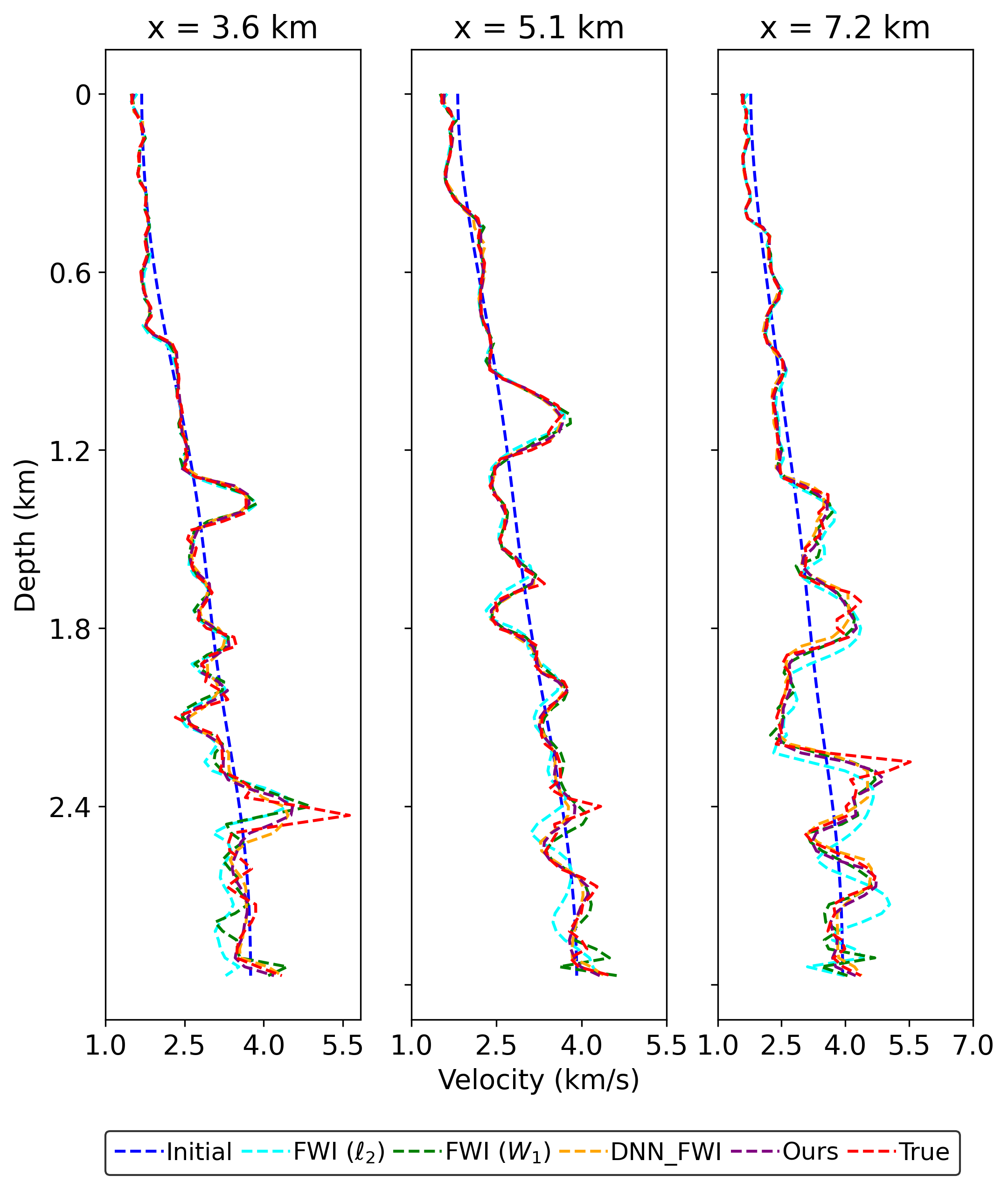}
		\end{minipage}
	}%
	\subfigure{
		\begin{minipage}[t]{0.5\linewidth}
			\centering
			\includegraphics[width=3in]{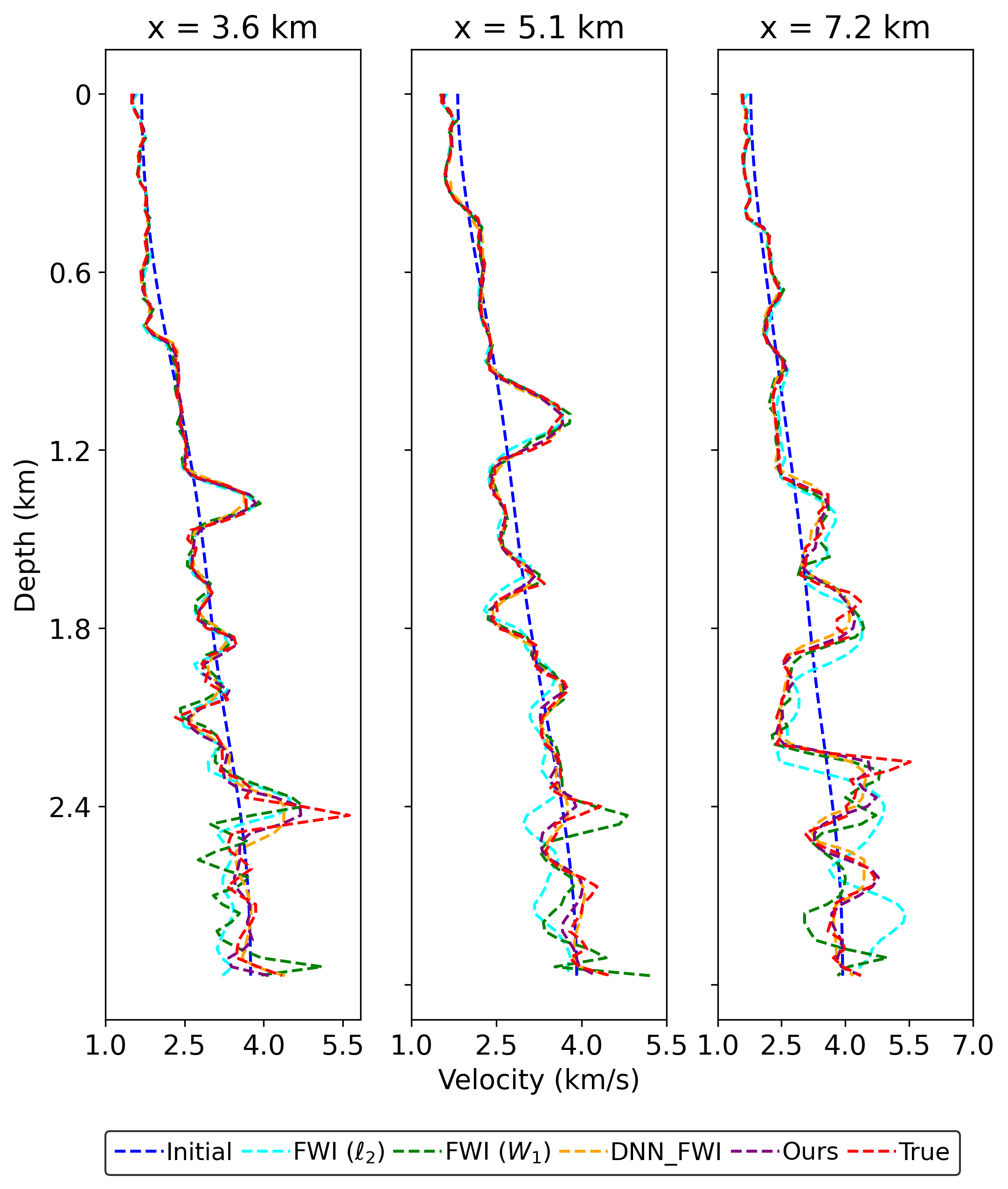}
		\end{minipage}
	}%
	\centering
	\caption{ The Marmousi velocity profiles of the true, initial, and inverted models
		obtained by different approaches at three horizontal locations. Left: noise-free observed data. Right: noisy observed data.}
	\label{fig_mar3}
\end{figure}

We further extract subsurface velocity model profiles at horizontal positions
of $x=3.6$ km, $x=5.1$ km, and $x=7.2$ km, as displayed in Figure \ref{fig_mar3}.
This comparative analysis underscores the superior performance of our method in
matching the true model. Quantitative metrics, including SNR, SSIM, and
relative $l_2$ error (as outlined in Table \ref{table_mar}), provide empirical
support for these findings.

\begin{figure}
	\centering
	\includegraphics[width=1.0\textwidth]{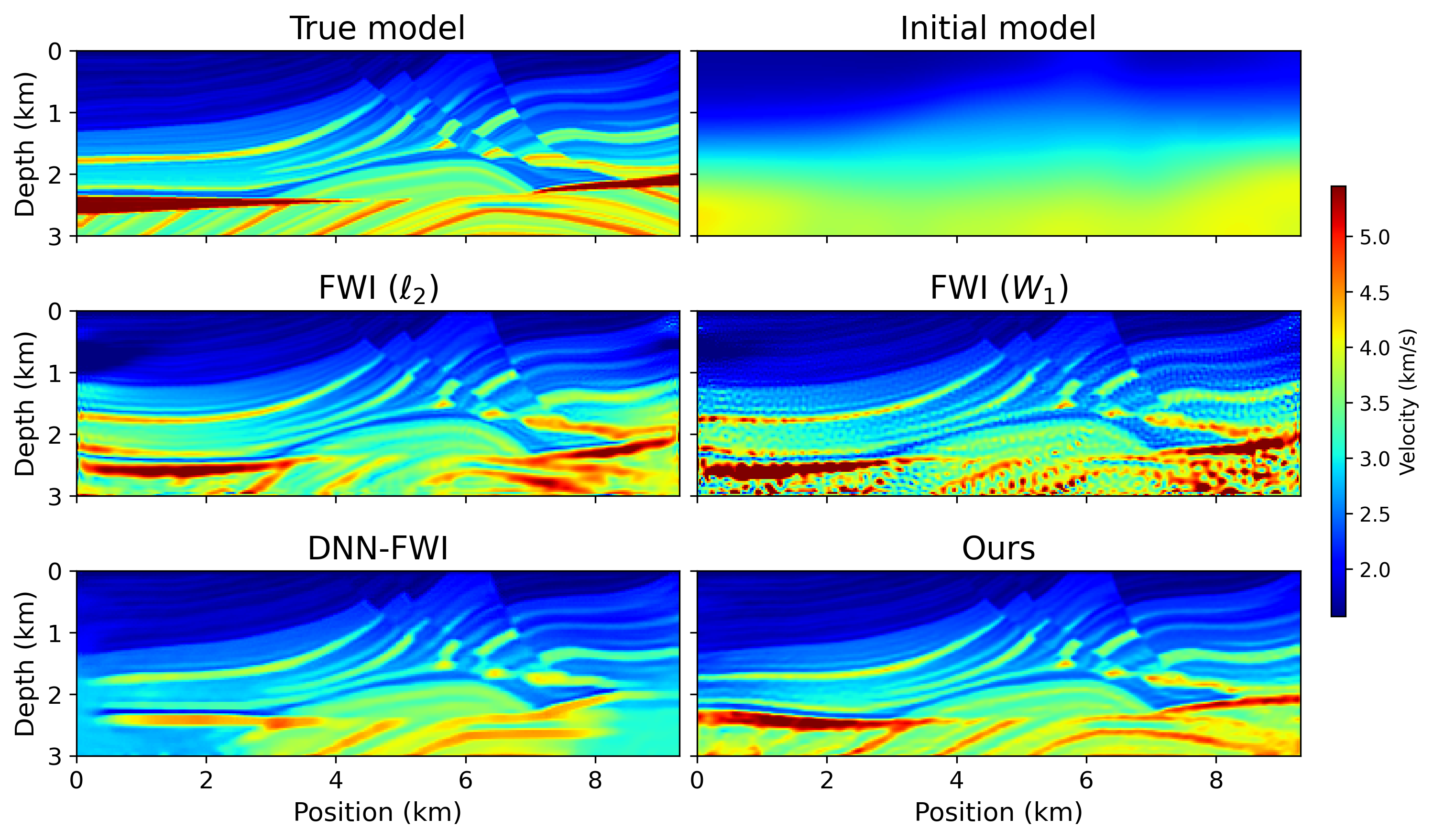}
	\caption{Inversion results for the noisy observed data of the Marmousi model.
		The top row displays the true velocity model and the initial guess.
		The second row shows the inverted results using traditional FWI constrained by the $l_2$-norm and the $1D-W_1$ distance.
		The third row presents the inverted results obtained through DNN-FWI and our method.}
	\label{fig_mar1}
\end{figure}

To showcase our method's capability in reconstructing the subsurface
velocity model with noisy measurement data, we conduct a comparison
of various inversion methods, as depicted in Figure \ref{fig_mar1}.
The results clearly indicate that, in comparison to the other three methods,
our approach yields superior numerical results for measurement data contaminated with noise.
This conclusion is further supported by the vertical profiles in Figure \ref{fig_mar3}
and the quantitative metrics presented in Table \ref{table_mar}.

\begin{figure}[H]
	\centering
	\includegraphics[width=1.0\textwidth]{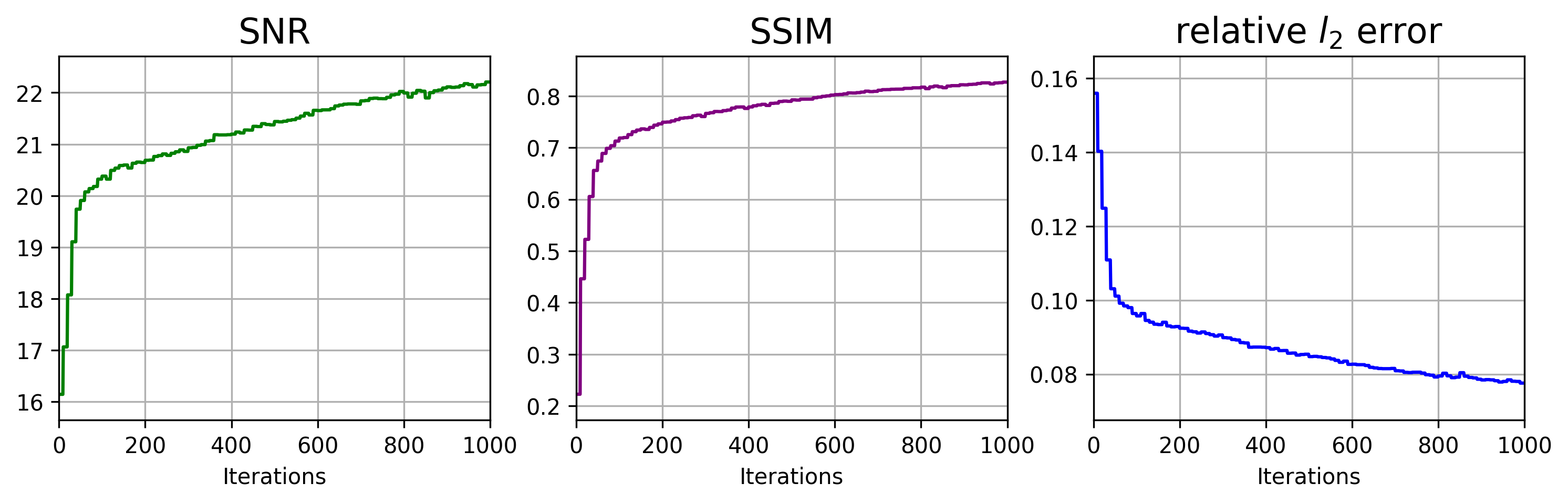}
	\caption{The trend of SNR, SSIM, and relative $l_2$ error of our method with
		respect to the number of iterations for the Marmousi model with noise-free observed data.}
	\label{fig_mar4}
\end{figure}

\begin{figure}[H]
	\centering
	\includegraphics[width=1.0\textwidth]{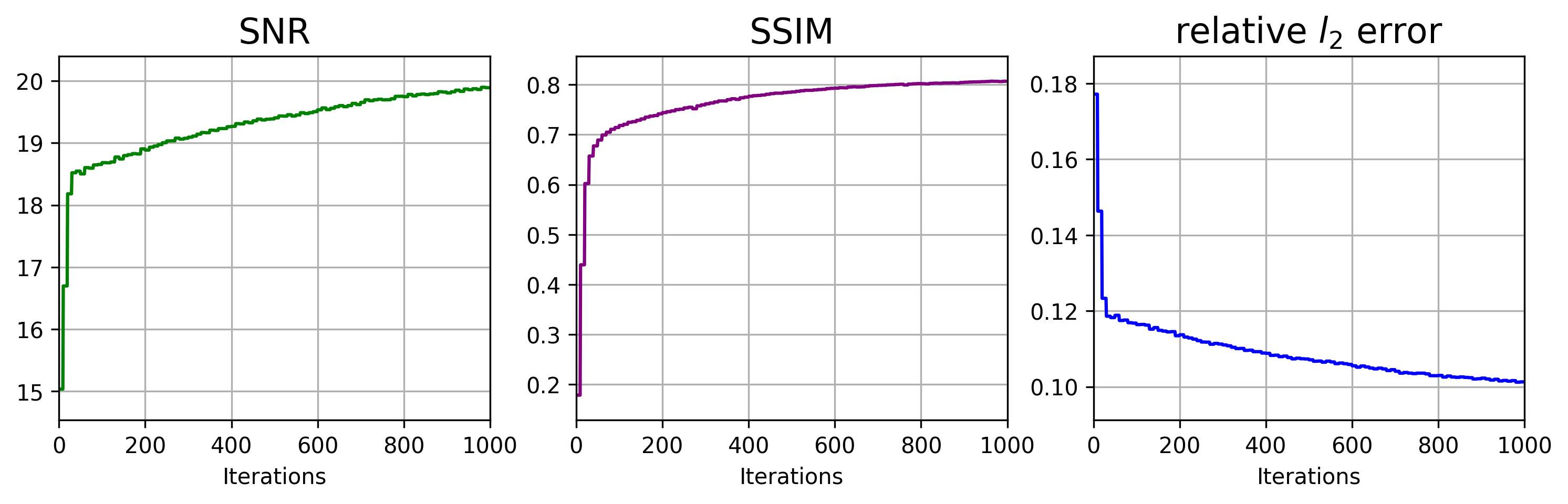}
	\caption{The trend of SNR, SSIM, and relative $l_2$ error of our method with
		respect to the number of iterations for the Marmousi model with noisy data.}
	\label{fig_mar5}
\end{figure}

Figures \ref{fig_mar4} and \ref{fig_mar5} depict how SNR, SSIM, and relative $l_2$ error evolve during the
optimization process, under both noisy and noise-free measurement conditions.
The consistent upward trends in SNR and SSIM, along with the consistent decrease
in relative $l_2$ error, signify the stability of our optimization process, which is a highly desirable
characteristic for practical applications. It's worth noting that the
optimization process displayed in Figures \ref{fig_mar4} and \ref{fig_mar5}
exhibits notably greater stability, suggesting a potentially more favorable
loss landscape for our method.

\begin{table}[H]
	\begin{center}
		\begin{spacing}{1.36}
			\resizebox{1\hsize}{!}{
				\begin{tabular}{cccccccccccccc}
					\hline
					\multirow{1}{*}{Case} & &\multirow{1}{*}{Index} & & \multirow{1}{*}{$\rm{FWI(l_2)}$} && \multicolumn{ 1}{c}{$\rm{FWI(W_1)}$} && \multicolumn{ 1}{c}{$\rm{DNN-FWI}$} && \multicolumn{1}{c}{Ours}&&  \\
					\cline{4-13}
					\hline
					&&SNR      &&15.41    &&17.81      &&16.86     &&\textbf{22.21}          \\
					Noise-free      &&SSIM     &&0.6118    &&0.7437     &&0.7900    &&\textbf{0.8266}       \\
					&&Error    &&0.1801   &&0.1287     &&0.1435    &&\textbf{0.0776}       \\
					
					\hline
					&&SNR    &&14.89         &&15.39      &&15.64       &&\textbf{19.89}         \\
					Noise         &&SSIM    &&0.5656        &&0.6084     &&0.7688      &&\textbf{0.8070}       \\
					&&Error    &&0.1801        &&0.1700     &&0.1653      &&\textbf{0.1012}     \\
					\hline
					
					\hline
				\end{tabular}
			}
			\caption{The SNR, SSIM, and relative $l_2$ error (Error) for the inverted Marmousi model
				obtained by different approaches.}\label{table_mar}
		\end{spacing}
	\end{center}
\end{table}
\subsection{Overthrust model}
The Overthrust model \cite{lecomte1994building} portrays a stratigraphy characterized by complex thrusting,
overlaying a previous extensional and rift sequence. The model exhibits varying
complexities, including a central thrust faulted anticline, an external monocline,
and a flat zone. The uppermost layer of the overthrust has undergone erosion and
is covered by a surface layer, symbolizing recent sediments.

\begin{figure}[H]
	\centering
	\includegraphics[width=1.0\textwidth]{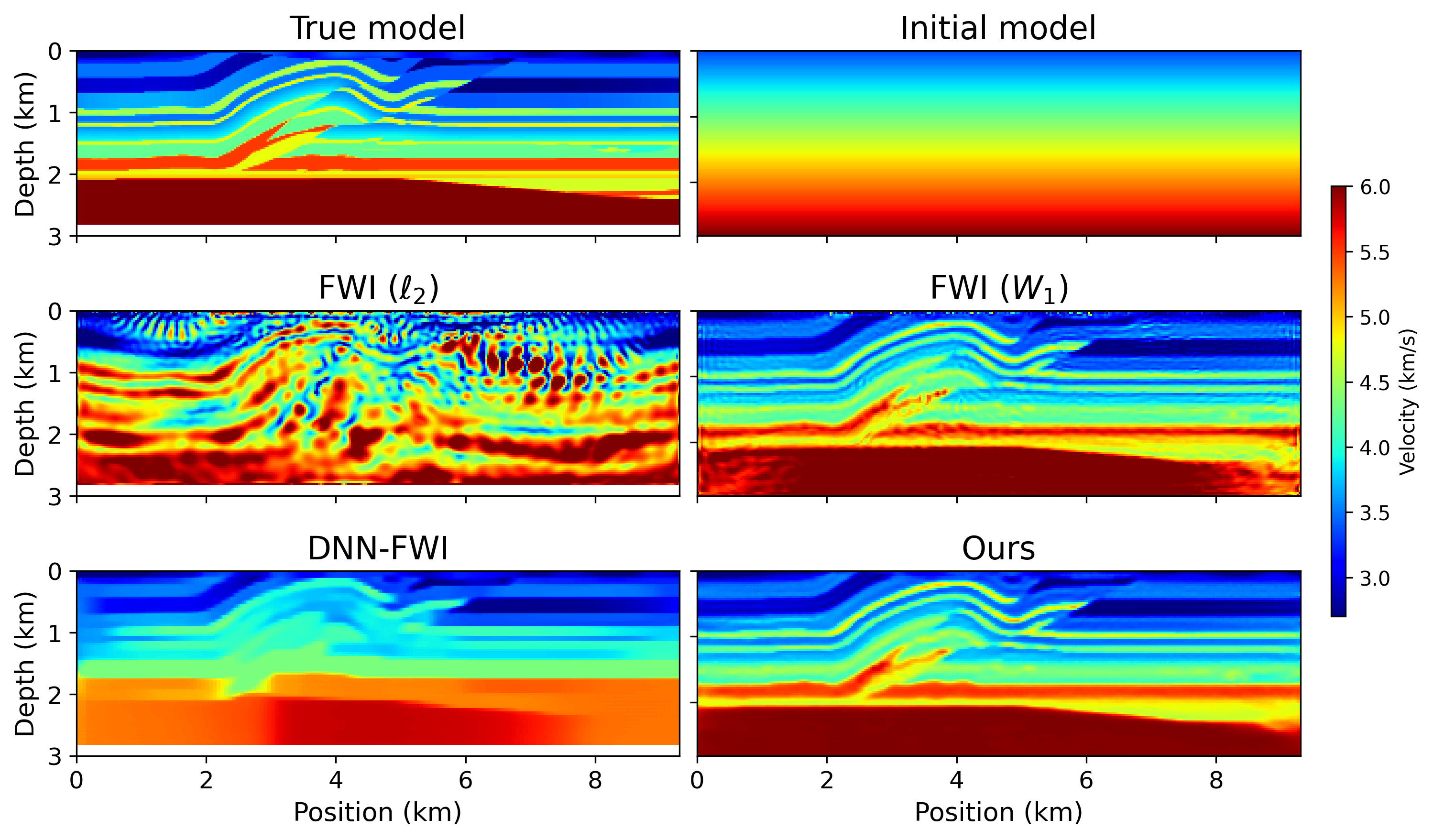}
	\caption{Inversion results for the noise-free observed data of the Overthrust model.
		The top row displays the true velocity model and the initial guess.
		The second row shows the inverted results using traditional FWI constrained by the $l_2$-norm and the $1D-W_1$ distance.
		The third row presents the inverted results obtained through DNN-FWI and our method.}
	\label{fig_over0}
\end{figure}

\begin{figure}[H]
	\centering
	\subfigure{
		\begin{minipage}[t]{0.5\linewidth}
			\centering
			\includegraphics[width=3in]{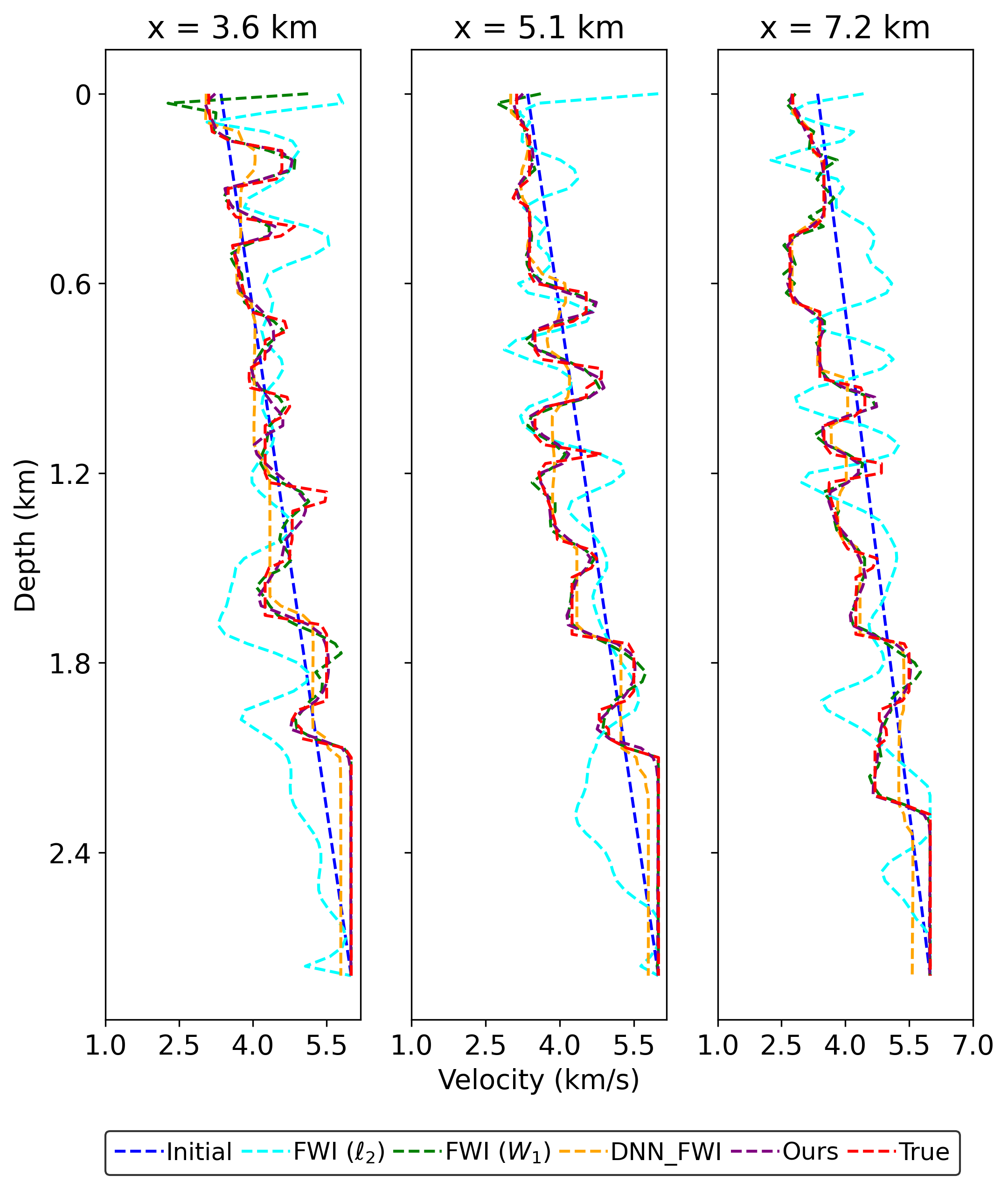}
		\end{minipage}
	}%
	\subfigure{
		\begin{minipage}[t]{0.5\linewidth}
			\centering
			\includegraphics[width=3in]{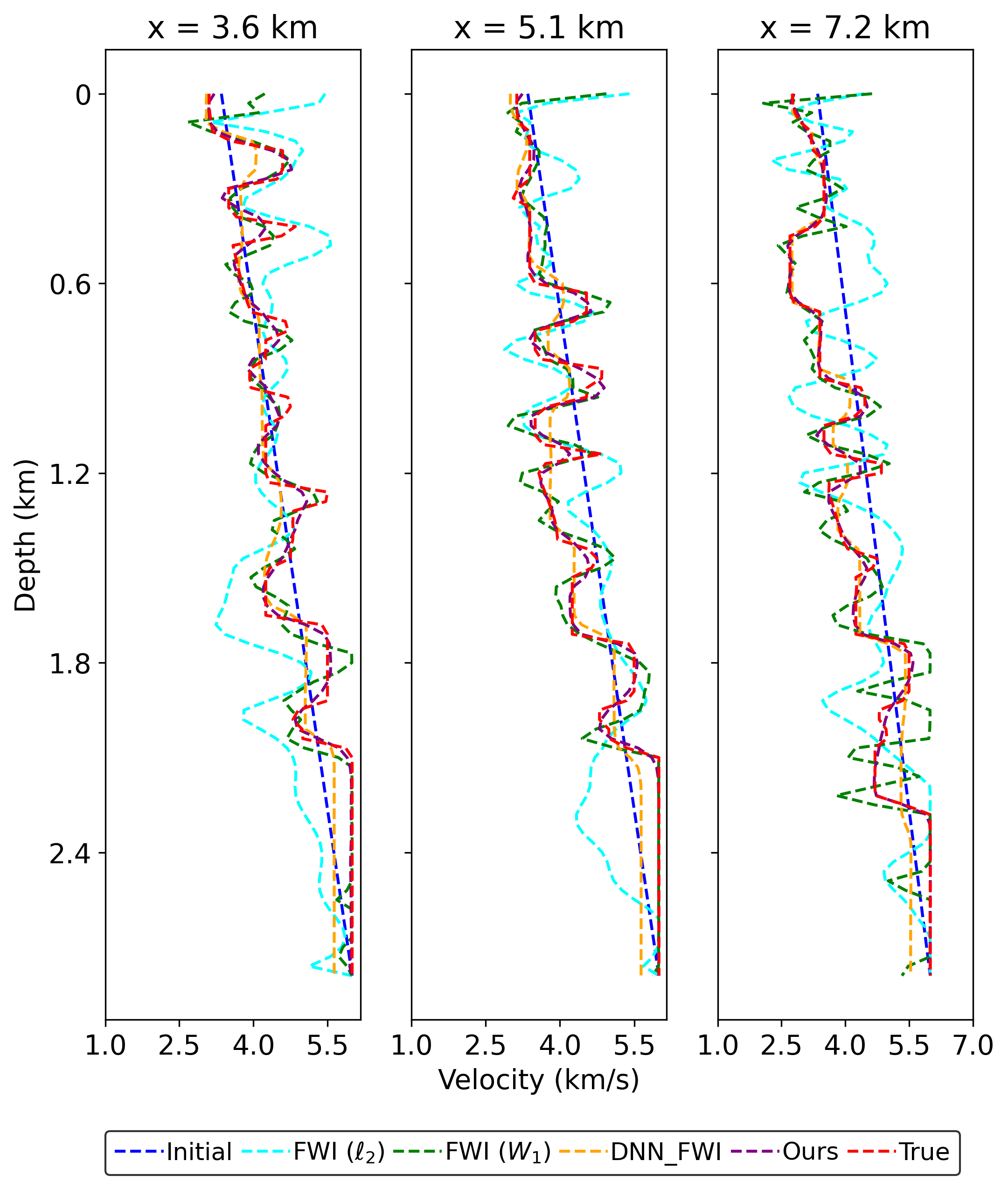}
		\end{minipage}
	}%
	\centering
	\caption{The Overthrust velocity profiles of the true, initial, and inverted models
		obtained by different approaches at three horizontal locations. Left: noise-free observed data. Right: noisy observed data.}\label{fig_over2}
\end{figure}

In Figure \ref{fig_over0}, we present the inversion results using noise-free measurement data.
It's evident that the FWI($l_2$) method struggles to reconstruct the subsurface velocity model.
In the case of DNN-FWI, only the upper region is accurately inverted.
In contrast, both FWI($W_1$) and our approach provide an improved solution,
accurately recovering the velocity in all regions. Likewise, when comparing
velocity profiles, as shown in Figure \ref{fig_over2}, the
velocity model obtained by our method aligns better with the true model compared to the others.

\begin{figure}[H]
	\centering
	\includegraphics[width=1.0\textwidth]{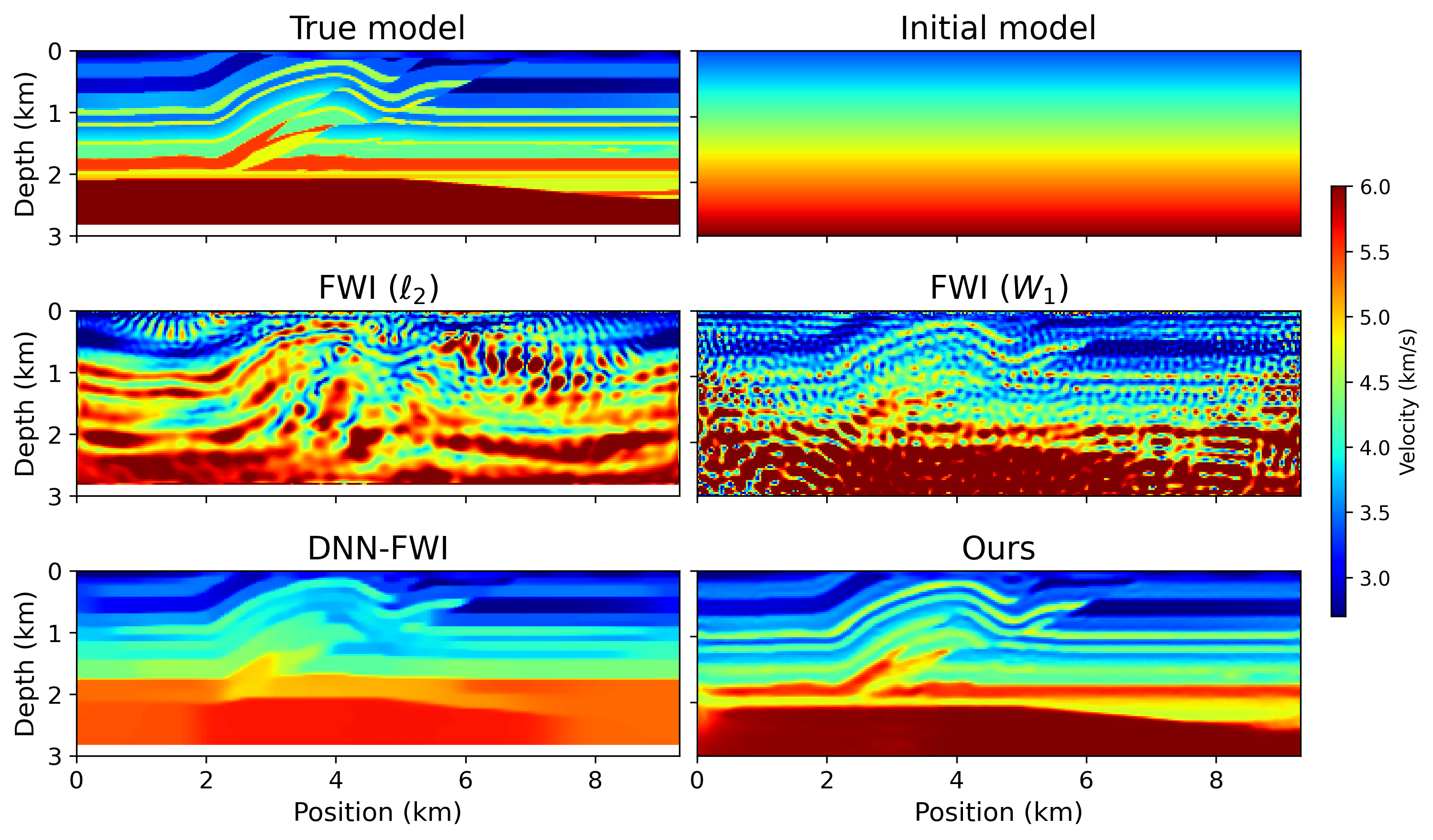}
	\caption{Inversion results for the noisy observed data of the Overthrust model.
		The top row displays the true velocity model and the initial guess.
		The second row shows the inverted results using traditional FWI constrained by the $l_2$-norm and the $1D-W_1$ distance.
		The third row presents the inverted results obtained through DNN-FWI and our method.}
	\label{fig_over1}
\end{figure}

Similarly, Figure \ref{fig_over1} provides a comparison of inversion results
using measurement data affected by random noise. It is evident that our method
outperforms all other compared approaches in terms of the inversion results.
The velocity profiles, as depicted in Figure \ref{fig_over2},
and the quantitative metrics listed in Table \ref{table_over},
once again emphasize the superiority of our method.

\begin{figure}[H]
	\centering
	\includegraphics[width=1.0\textwidth]{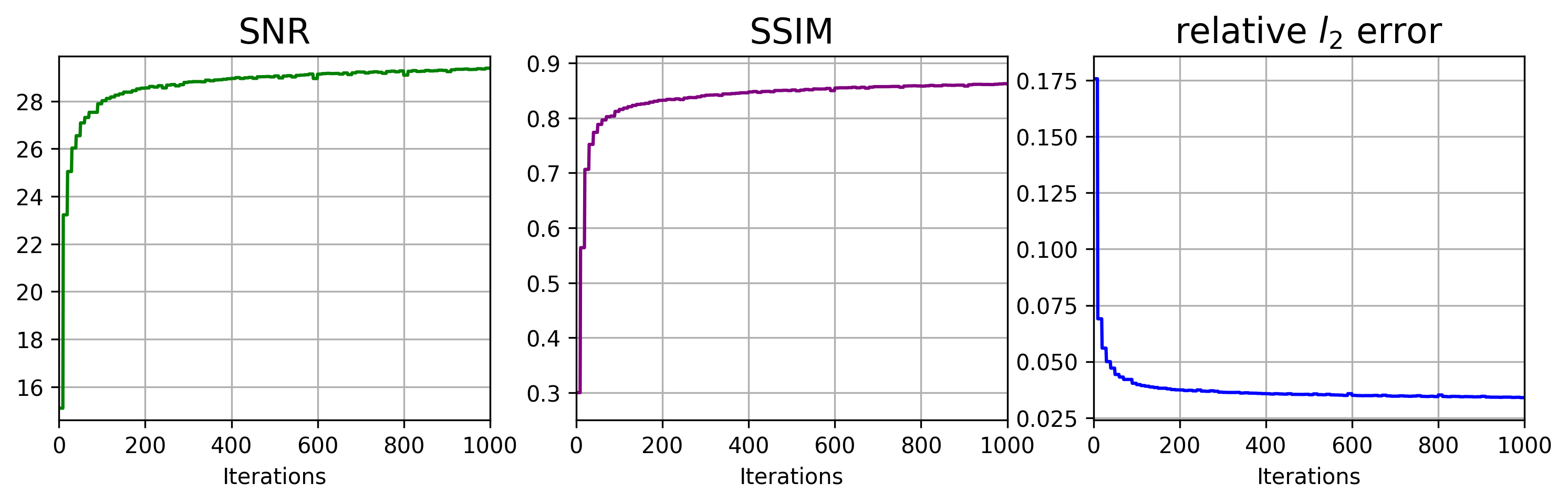}
	\caption{The trend of SNR, SSIM, and relative $l_2$ error of our method with
		respect to the number of iterations for the Overthrust model with noise-free observed data.}
	\label{fig_over3}
\end{figure}

\begin{figure}[H]
	\centering
	\includegraphics[width=1.0\textwidth]{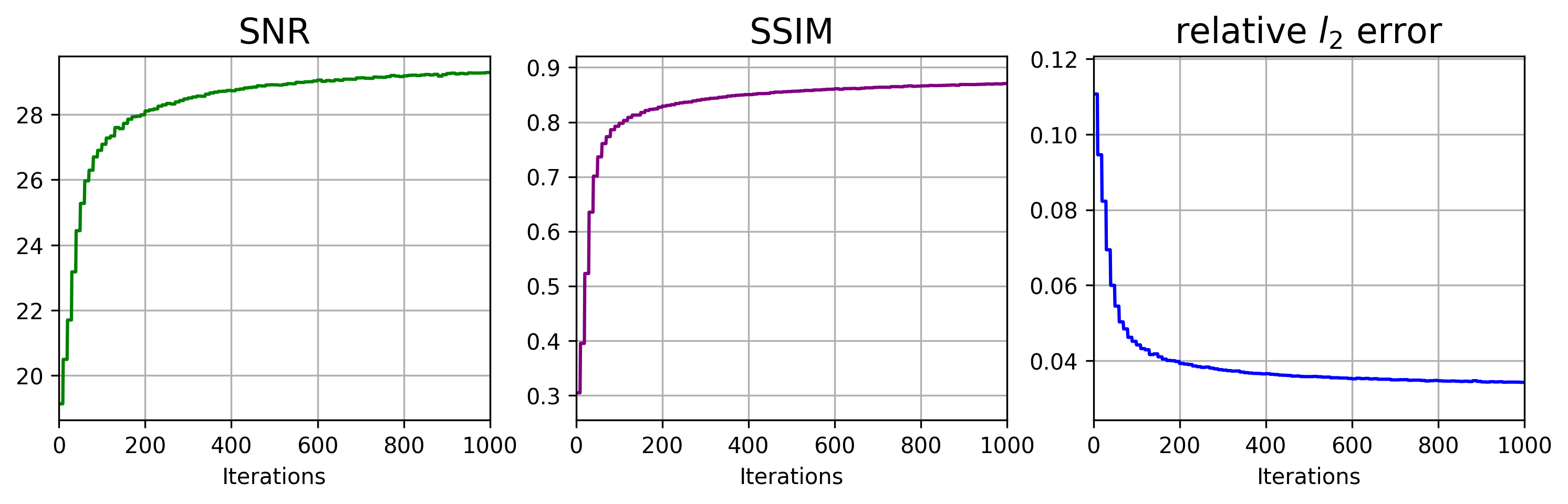}
	\caption{The trend of SNR, SSIM, and relative $l_2$ error of our method with
		respect to the number of iterations for the Overthrust model with noisy observed data.}
	\label{fig_over4}
\end{figure}

Furthermore, we have presented convergence curves in
Figures \ref{fig_over3} and \ref{fig_over4} that track the progress of SNR, SSIM,
and relative $l_2$ error throughout the optimization process facilitated by our method.
These curves vividly illustrate the rapid convergence achieved by our inversion approach.
Moreover, the elevated values of SNR and SSIM, coupled with the low relative $l_2$ error,
underscore the high precision of the inverted velocity model obtained by our method.

\begin{table}[H]
	\begin{center}
		\begin{spacing}{1.36}
			
			\resizebox{1\hsize}{!}{
				\begin{tabular}{cccccccccccccc}
					\hline
					\multirow{1}{*}{Case} & &\multirow{1}{*}{Index} & & \multirow{1}{*}{$\rm{FWI(l_2)}$} && \multicolumn{ 1}{c}{$\rm{FWI(W_1)}$} && \multicolumn{ 1}{c}{$\rm{DNN-FWI}$} && \multicolumn{1}{c}{Ours}&&  \\
					\cline{4-13}
					\hline
					&&SNR      &&14.10     &&27.43      &&22.24     &&\textbf{29.39}          \\
					Noise-free      &&SSIM     &&0.0792    &&0.8332     &&0.6266    &&\textbf{0.8624}       \\
					&&Error    &&0.1973    &&0.0422     &&0.0773    &&\textbf{0.0339}       \\
					
					\hline
					&&SNR    &&14.12         &&18.22      &&21.77       &&\textbf{29.29}         \\
					Noise         &&SSIM    &&0.0810        &&0.4420     &&0.6023      &&\textbf{0.8705}       \\
					&&Error    &&0.1968        &&0.1228     &&0.0815      &&\textbf{0.0343}     \\
					\hline
					
					\hline
				\end{tabular}
			}
			\caption{The SNR, SSIM, and relative $l_2$ error (Error) for the inverted Overthrust model
				obtained by different approaches.}\label{table_over}
		\end{spacing}
	\end{center}
\end{table}
\subsection{Marmousi2 model}
\begin{figure}[H]
	\centering
	\includegraphics[width=1.0\textwidth]{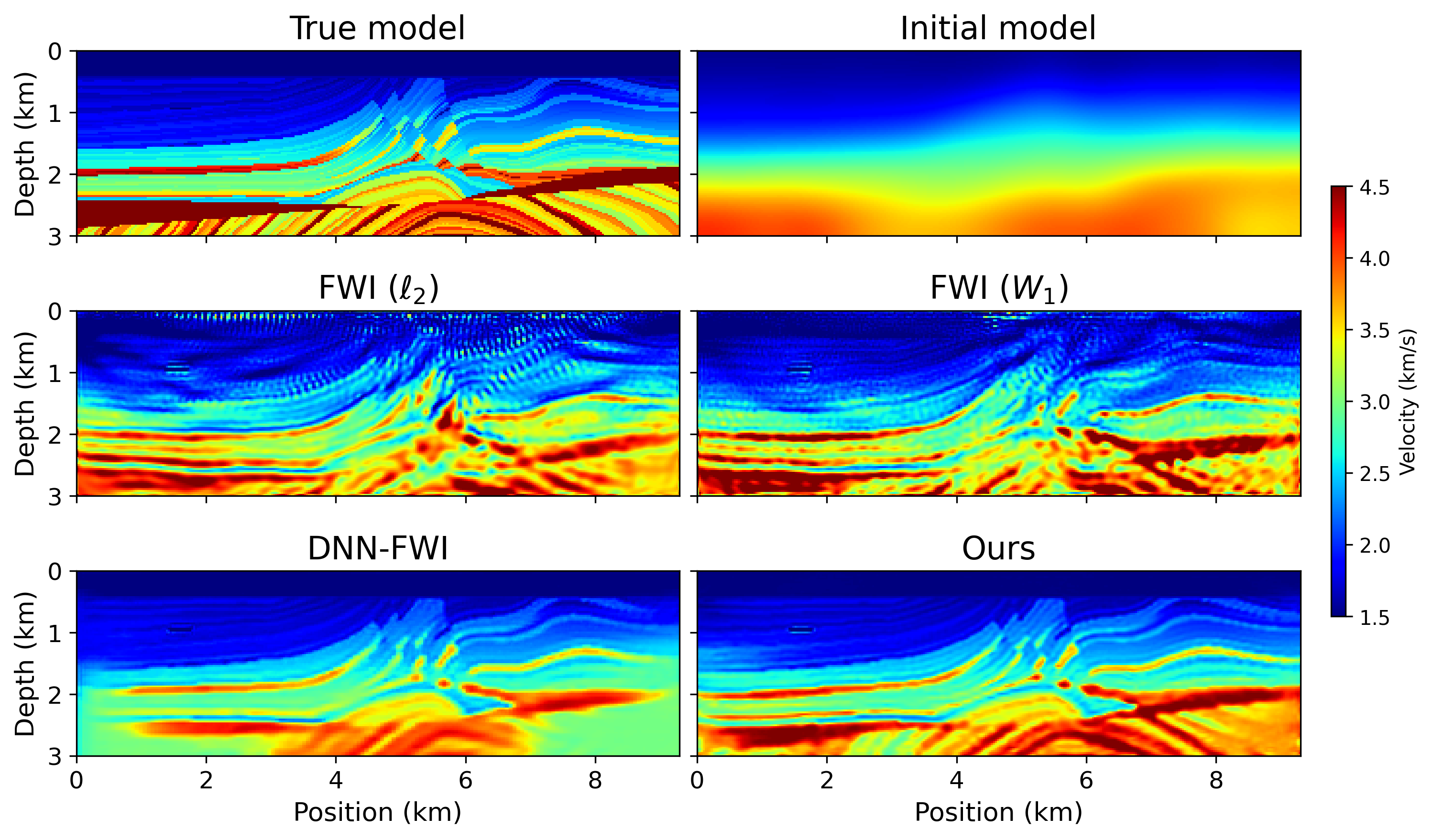}
	\caption{Inversion results for the noise-free observed data of the Marmousi2 model.
		The top row displays the true velocity model and the initial guess.
		The second row shows the inverted results using traditional FWI constrained by the $l_2$-norm and the $1D-W_1$ distance.
		The third row presents the inverted results obtained through DNN-FWI and our method.}
	\label{fig_mar2-0}
\end{figure}

Lastly, we further assess the effectiveness of our inversion method using the Marmousi2 model,
an extension of the original Marmousi model. The Marmousi2 model encompasses a broader region,
with the original Marmousi model positioned near its center. This extended model incorporates
a range of structurally simple but stratigraphically complex features,
surpassing the complexity of the original version. The downsampled version of the Marmousi2 used
in our experiments has dimensions of $(z\times x)=(100\times 300)$.
The velocity in the model ranges from 1140 m/s to 4700 m/s, and the true model is depicted in Figure \ref{fig_mar2-0}.

\begin{figure}[H]
	\centering
	\includegraphics[width=1.0\textwidth]{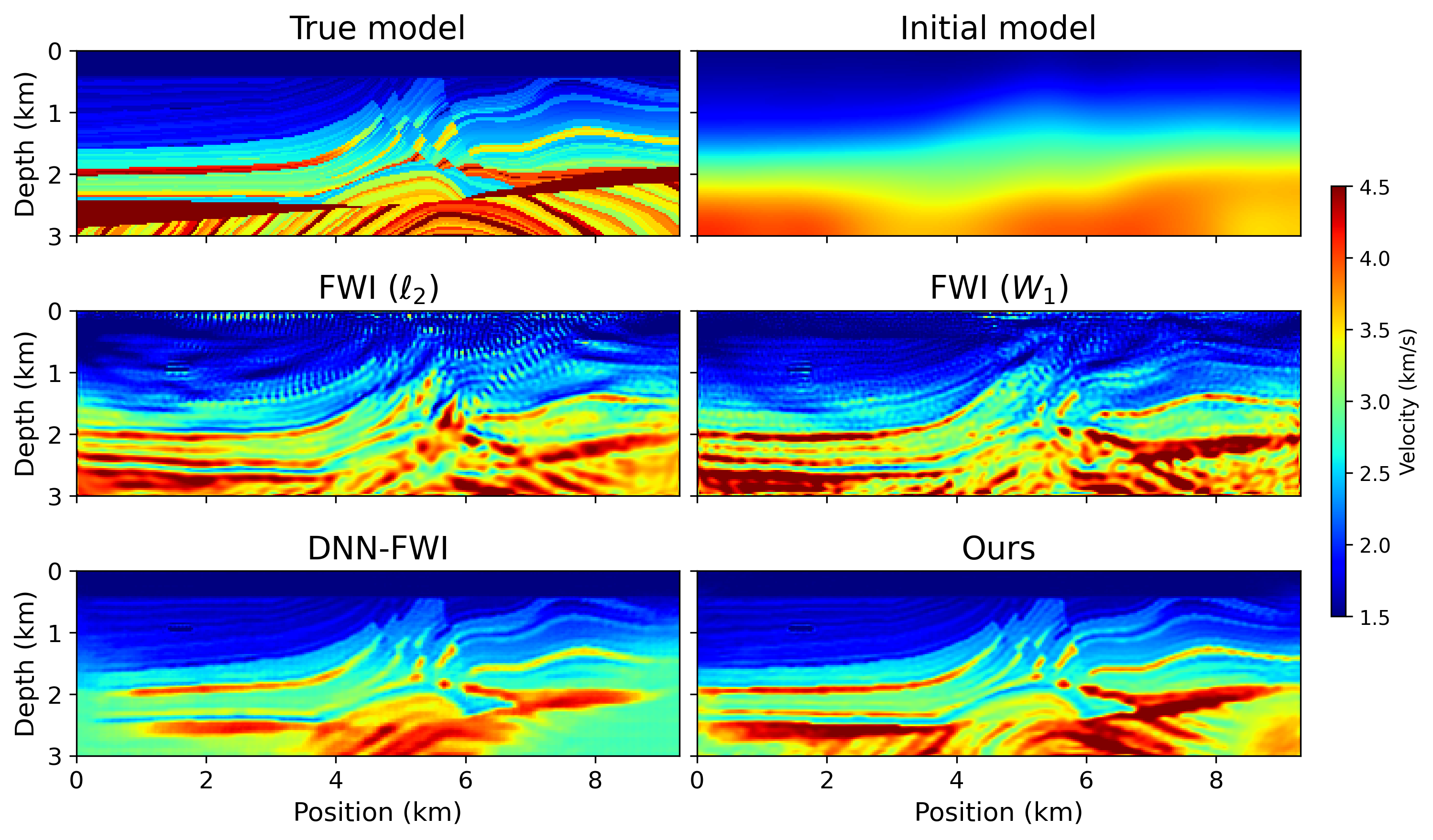}
	\caption{Inversion results for the noisy observed data of the Marmousi2 model.
		The top row displays the true velocity model and the initial guess.
		The second row shows the inverted results using traditional FWI constrained by the $l_2$-norm and the $1D-W_1$ distance.
		The third row presents the inverted results obtained through DNN-FWI and our method.}
	\label{fig_mar2-1}
\end{figure}

Figure \ref{fig_mar2-0} displays the velocity results obtained by different inversion methods.
Traditional FWI($l_2$) and FWI($W_1$) methods offer only approximate velocity estimates
and struggle to accurately reconstruct the detail of the true model. DNN-FWI improves inversion
performance, accurately reconstructing the upper part of the velocity model,
but faces challenges in delivering accurate results for deeper regions.
Our method, on the other hand, consistently provides more precise velocity and
structural representations for the Marmousi2 model, especially in the deeper
layers. These findings hold true even when working with noisy measurement
data, as illustrated in Figure \ref{fig_mar2-1}.

\begin{figure}[H]
	\centering
	\subfigure{
		\begin{minipage}[t]{0.5\linewidth}
			\centering
			\includegraphics[width=3in]{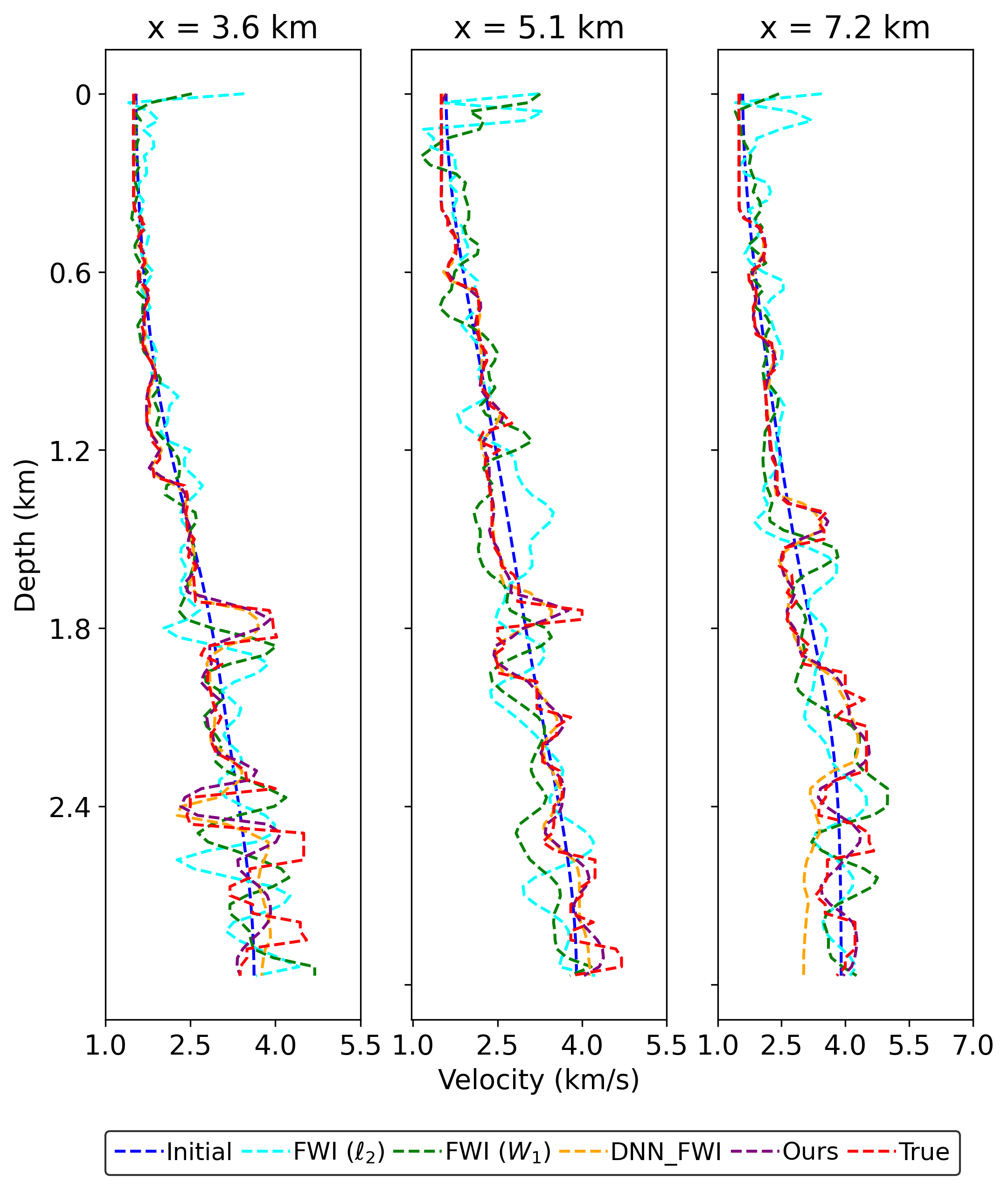}
		\end{minipage}
	}%
	\subfigure{
		\begin{minipage}[t]{0.5\linewidth}
			\centering
			\includegraphics[width=3in]{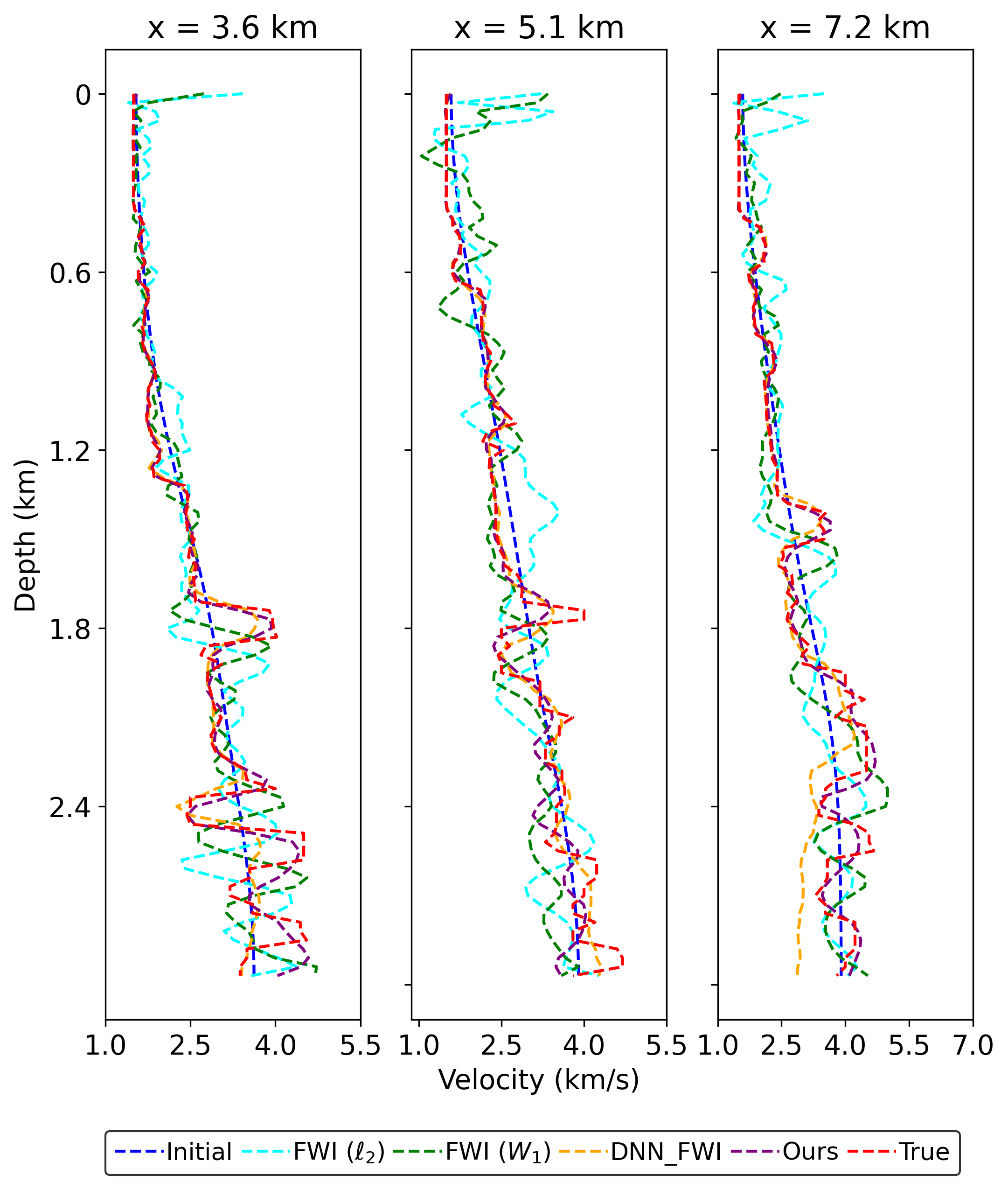}
		\end{minipage}
	}%
	\centering
	\caption{The Marmousi2 velocity profiles of the true, initial, and inverted models
		obtained by different approaches at three horizontal locations.
		Left: noise-free observed data. Right: noisy observed data.}\label{fig_mar2-2}
\end{figure}

Figure \ref{fig_mar2-2} illustrates a comparison of vertical velocity profiles
obtained by different inversion methods. Quantitative metrics, including SNR, SSIM, and
relative $l_2$ error, are summarized in Table \ref{table_mar2}.
These observations demonstrate the consistent and stable performance of
our method in addressing the full-waveform inversion problem.
Importantly, we observe that our proposed approach delivers
consistent accurate inversion results in scenarios with and without noise,
affirming its robustness to noise.

\begin{figure}[H]
	\centering
	\includegraphics[width=1.0\textwidth]{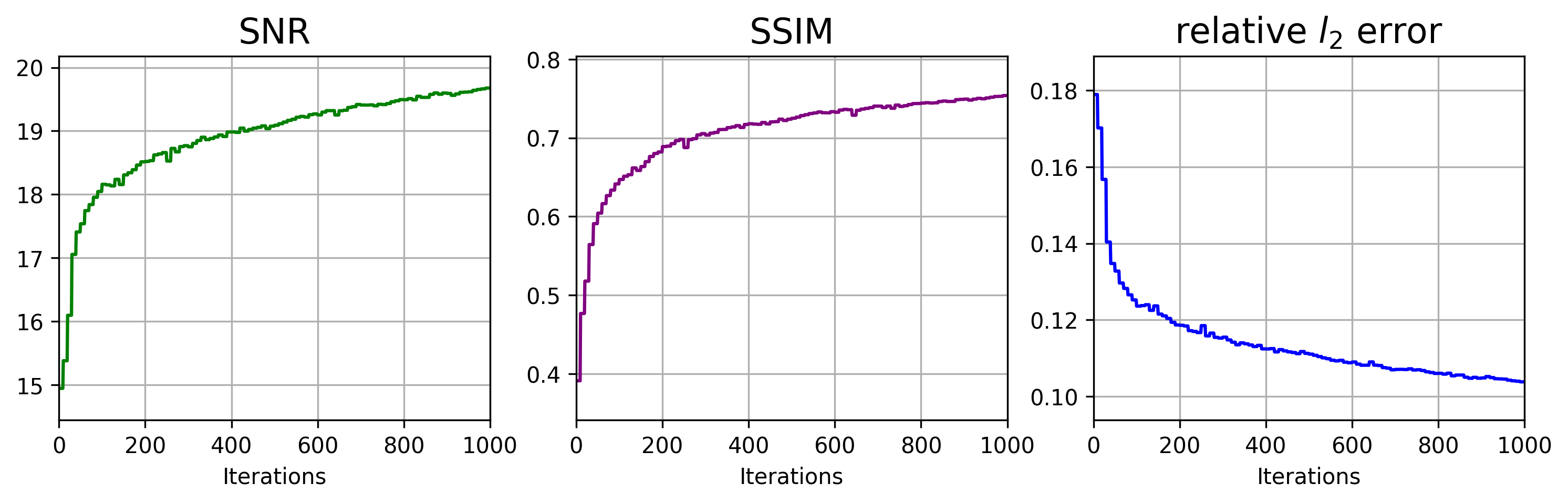}
	\caption{The trend of SNR, SSIM, and relative $l_2$ error of our method with
		respect to the number of iterations for the Marmousi2 model with noise-free observed data.}
	\label{fig_mar2-3}
\end{figure}

\begin{figure}[H]
	\centering
	\includegraphics[width=1.0\textwidth]{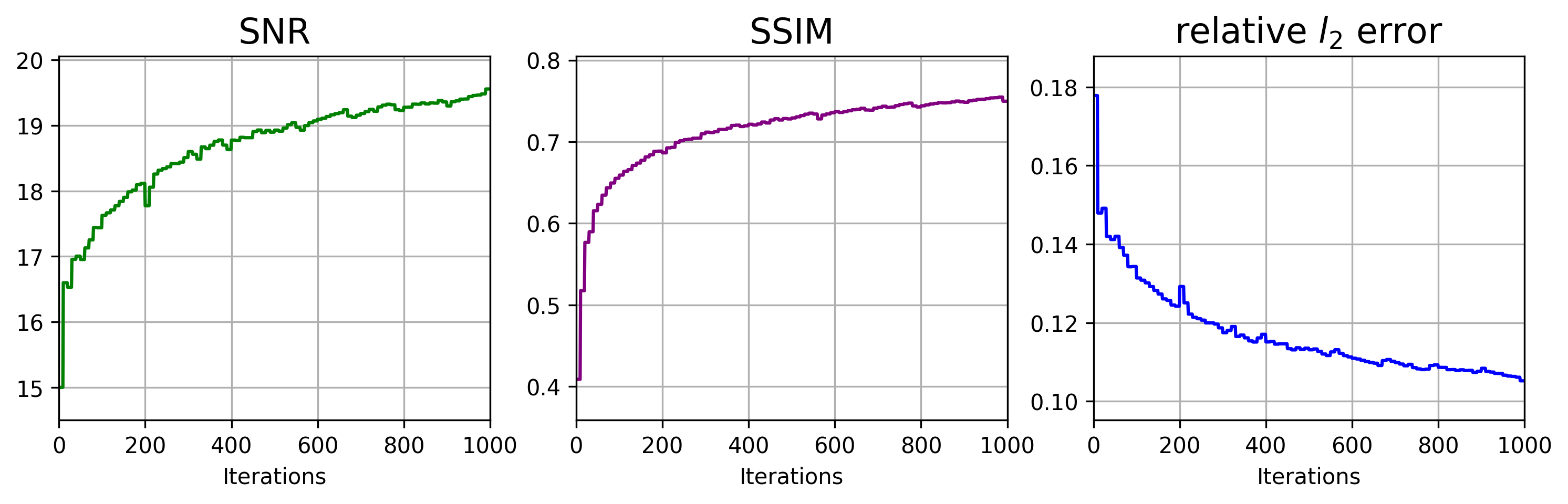}
	\caption{The trend of SNR, SSIM, and relative $l_2$ error of our method with
		respect to the number of iterations for the Marmousi2 model with noisy observed data.}
	\label{fig_mar2-4}
\end{figure}

Similarly, we have also plotted the convergence curves for SNR, SSIM, and relative $l_2$
error of the iterations progress of our method, as shown in Figures \ref{fig_mar2-3}, \ref{fig_mar2-4}.
These stable trends further demonstrate the reliability of the inversion method provided by us.


\begin{table}[H]
	\begin{center}
		\begin{spacing}{1.36}
			
			\resizebox{1\hsize}{!}{
				\begin{tabular}{cccccccccccccc}
					\hline
					\multirow{1}{*}{Case} & &\multirow{1}{*}{Index} & & \multirow{1}{*}{$\rm{FWI(l_2)}$} && \multicolumn{ 1}{c}{$\rm{FWI(W_1)}$} && \multicolumn{ 1}{c}{$\rm{DNN-FWI}$} && \multicolumn{1}{c}{Ours}&&  \\
					\cline{4-13}
					\hline
					&&SNR      &&13.25     &&13.67      &&17.14     &&\textbf{19.68}          \\
					Noise-free      &&SSIM     &&0.1270    &&0.2003     &&0.7261    &&\textbf{0.7537}       \\
					&&Error    &&0.2176    &&0.2073     &&0.1390    &&\textbf{0.1038}       \\
					
					\hline
					&&SNR     &&13.27         &&13.59      &&16.1317       &&\textbf{19.56}         \\
					Noise         &&SSIM    &&0.1257        &&0.1870     &&0.7069      &&\textbf{0.7496}       \\
					&&Error     &&0.2170        &&0.2091     &&0.1561      &&\textbf{0.1052}     \\
					\hline
					
					\hline
				\end{tabular}
			}
			\caption{The SNR, SSIM, and relative $l_2$ error (Error) for the inverted Marmousi2 model
				obtained by different approaches.}\label{table_mar2}
		\end{spacing}
	\end{center}
\end{table}

\section{Conclusion}\label{sec5}
In this paper, we propose an unsupervised full-waveform inversion method
that utilizes re-parametrization through a neural network with random weights to approximate the posterior distribution of the Bayesian inverse problem. This approach inherits the benefits of unsupervised
learning, specifically, it does not necessitate external training samples,
offering flexibility and practical ease of use.

We validate our approach using the established Marmousi, Marmousi2, and Overthrust models,
both with noise-free and noisy observed data. Our experiments demonstrate
superior performance compared to two representative non-learning-based methods
and one recent unsupervised deep-learning method.
Significantly, we have observed that the provided inversion method exhibits robustness to measurement noise. In other words, even when the measurement data contains noise, our inversion method consistently produces accurate results.

Several promising avenues for further research are available.
Firstly, our inversion method relies on an initial guess.
We can explore techniques like model-agnostic meta-learning
and adversarial pretraining to improve initialization.
Secondly, our method has so far been tested on benchmark
models. Our future plans involve applying it to practical model inversion problems.

\bigskip
\noindent{\bf Acknowledgments}\\
This work is sponsored by the National Key R\&D Program of China  Grant No. 2022YFA1008200 (Z. X.) and No. 2020YFA0712000 (Z. M.), the Shanghai Sailing Program (Z. X.), the Natural Science Foundation of Shanghai Grant No. 20ZR1429000  (Z. X.), the National Natural Science Foundation of China Grant No. 62002221 (Z. X.), the National Natural Science Foundation of China Grant No. 12101401 (Z. M.), the National Natural Science Foundation of China Grant No. 12031013 (Z. M.), Shanghai Municipal of Science and Technology Major Project No. 2021SHZDZX0102, and the HPC of School of Mathematical Sciences and the Student Innovation Center, and the Siyuan-1 cluster supported by the Center for High Performance Computing at Shanghai Jiao Tong University.
\medskip

\bibliographystyle{unsrt}
\bibliography{dip_wave_v2}

\end{document}